\documentclass[11pt,a4paper]{article}
\usepackage[cp1251]{inputenc}
\usepackage{amsfonts}
\usepackage{amssymb}
\usepackage{amsthm}
\usepackage{amsmath}
\usepackage{mathtext}
\usepackage{longtable}
\usepackage[russian]{babel}
\usepackage{color}

\title{Multivariate scale-mixed stable
distributions and related limit theorems\thanks{Supported by Russian
Science Foundation, project 18-11-00155.}}
\author{Victor Korolev\thanks{Faculty of Computational Mathematics
and Cybernetics, Moscow State University, Moscow, Russia; Hanghzhou
Dianzi University, Hangzhou, China; Institute of Informatics
Problems, Federal Research Center ``Computer Science and Control'',
Russian Academy of Sciences, Moscow, Russia. E-mail:
vkorolev@cs.msu.ru}, Alexander Zeifman\thanks{Vologda State
University, Vologda, Russia; Institute of Informatics Problems,
Federal Research Center ``Computer Science and Control'',  Russian
Academy of Sciences, Moscow, Russia; Vologda Research Center of the
Russian Academy of Sciences, Vologda, Russia. E-mail:
a$\_$zeifman@mail.ru}}
\date{}

\renewcommand{\r}{\mathbb{R}}

\newcommand{\pto}{\stackrel{P}{\longrightarrow}}
\newcommand{\eqd}{\stackrel{d}{=}}

\renewcommand{\le}{\leqslant}
\renewcommand{\ge}{\geqslant}
\newcommand{\nb}{N\!\!\!B}

\begin{document}

\sloppy

\maketitle

%\renewcommand{\contentsname}{Contents}

%\tableofcontents

\begin{center} {\bf Abstract} \end{center}

In the paper, multivariate probability distributions are considered
that are representable as scale mixtures of multivariate
elliptically contoured stable distributions. It is demonstrated that
these distributions form a special subclass of scale mixtures of
multivariate elliptically contoured normal distributions. Some
properties of these distributions are discussed. Main attention is
paid to the representations of the corresponding random vectors as
products of independent random variables. In these products,
relations are traced of the distributions of the involved terms with
popular probability distributions. As examples of distributions of
the class of scale mixtures of multivariate elliptically contoured
stable distributions, multivariate generalized Linnik distributions
are considered in detail. Their relations with multivariate
`ordinary' Linnik distributions, multivariate normal, stable and
Laplace laws as well as with univariate Mittag-Leffler and
generalized Mittag-Leffler distributions are discussed. Limit
theorems are proved presenting necessary and sufficient conditions
for the convergence of the distributions of random sequences with
independent random indices (including sums of a random number of
random vectors and multivariate statistics constructed from samples
with random sizes) to scale mixtures of multivariate elliptically
contoured stable distributions. The property of scale-mixed
multivariate stable distributions to be both scale mixtures of a
non-trivial multivariate stable distribution and a normal scale
mixture is used to obtain necessary and sufficient conditions for
the convergence of the distributions of random sums of random
vectors with both infinite or finite covariance matrices to the
multivariate generalized Linnik distribution.
\\[5pt]

%$\emph{\textbf{x}}$

\section{Introduction}

In the paper, multivariate probability distributions are considered
that are representable as scale mixtures of multivariate
elliptically contoured stable distributions. It is demonstrated that
each of these distributions can be represented as a scale mixture of
multivariate elliptically contoured normal distributions. On the
other hand, since the multivariate normal distribution is stable
with $\alpha=2$, any multivariate normal scale mixture can be
regarded as a `trivial' multivariate scale-mixed stable
distribution. Most results presented in the paper concern scale
mixtures of `non-trivial' multivariate stable laws with
$0<\alpha<2$. Some properties of these mixtures are discussed. Main
attention is paid to the representations of the corresponding random
vectors as products of independent random variables. In these
products, relations are traced of the distributions of the involved
terms with popular probability distributions. As examples of
distributions of the class of scale mixtures of multivariate
elliptically contoured stable distributions, multivariate
generalized Linnik distributions are considered in detail. Limit
theorems are proved presenting necessary and sufficient conditions
for the convergence of the distributions of random sequences with
independent random indices (including sums of a random number of
random vectors and multivariate statistics constructed from samples
with random sizes) to scale mixtures of multivariate elliptically
contoured stable distributions. As particular cases, conditions are
obtained for the convergence of the distributions of random sums of
random vectors with both infinite or finite covariance matrices to
the multivariate generalized Linnik distribution.

Along with general properties of the class of scale mixtures of
multivariate elliptically contoured stable distributions, some
important and popular special cases are considered in detail. We
study the multivariate (generalized) Linnik and related
(generalized) Mittag-Leffler distributions, their interrelation and
their relations with multivariate `ordinary' Linnik distributions,
multivariate normal, stable and Laplace laws as well as with
univariate `ordinary' Mittag-Leffler distributions. Namely, we
consider mixture representations for the generalized Mittag-Leffler
and multivariate generalized Linnik distributions. We continue the
research we started in \cite{KorolevZeifman2016, KorolevZeifmanKMJ,
Korolevetal2018, Korolevetal2019}. In most papers (see, e. g.,
\cite{Anderson1993, Devroye1990, KotzOstrovskiiHayfavi1995a,
KotzOstrovskiiHayfavi1995b, KotzOstrovskii1996, Kotz2001,
Kozubowski1998, Kozubowski1999, Lin1998, LimTeo2009, Ostrovskii1995,
Pakes1998, Pillai1985}), the properties of the (multivariate)
generalized Linnik distribution and of the Mittag-Leffler
distributions were deduced by analytical methods from the properties
of the corresponding probability densities and/or characteristic
functions. Instead, here we use the approach which can be regarded
as arithmetical in the space of random variables or vectors. Within
this approach, instead of the operation of scale mixing in the space
of distributions, we consider the operation of multiplication in the
space of random vectors/variables provided the multipliers are
independent. This approach considerably simplifies the reasoning and
makes it possible to notice some general features of the
distributions under consideration. We prove mixture representations
for general scale mixtures of multivariate elliptically contoured
stable distributions and their particular cases in terms of normal,
Laplace, generalized gamma (including exponential, gamma and
Weibull) and stable laws and establish the relationship between the
mixing distributions in these representations. In particular, we
prove that the multivariate generalized Linnik distribution is a
multivariate normal scale mixture with the generalized
Mittag-Leffler mixing distribution and, moreover, this
representation can be used as the definition of the multivariate
generalized Linnik distribution. Based on these representations, we
prove some limit theorems for random sums of independent random
vectors with both infinite and finite covariance matrices. As a
particular case, we prove some theorems in which the multivariate
generalized Linnik distribution plays the role of the limit law. By
doing so, we demonstrate that the scheme of geometric (or, in
general, negative binomial) summation is far not the only asymptotic
setting (even for sums of independent random variables) in which the
multivariate generalized Linnik law appears as the limit
distribution.

In \cite{KorolevZeifmanKMJ} we showed that along with the
traditional and well-known representation of the univariate Linnik
distribution as the scale mixture of a strictly stable law with
exponential mixing distribution, there exists another representation
of the Linnik law as the normal scale mixture with the
Mittag-Leffler mixing distribution. The former representation makes
it possible to treat the Linnik law as the limit distribution for
geometric random sums of independent identically distributed random
variables (random variables) in which summands have infinite
variances. The latter normal scale mixture representation opens the
way to treating the Linnik distribution as the limit distribution in
the central limit theorem for random sums of independent random
variables in which summands have {\it finite} variances. Moreover,
being scale mixtures of normal laws, the Linnik distributions can
serve as the one-dimensional distributions of a special subordinated
Wiener process often used as models of the evolution of stock prices
and financial indexes. Strange as it may seem, the results
concerning the possibility of representation of the Linnik
distribution as a scale mixture of normals were never explicitly
presented in the literature in full detail before the paper
\cite{KorolevZeifmanKMJ} saw the light, although the property of the
Linnik distribution to be a normal scale mixture is something almost
obvious. Perhaps, the paper \cite{KotzOstrovskii1996} was the
closest to this conclusion and exposed the representability of the
Linnik law as a scale mixture of Laplace distributions with the
mixing distribution written out explicitly. These results became the
base for our efforts to extend them from the Linnik distribution to
the multivariate generalized Linnik law and more general scale
mixtures of multivariate stable distributions. Methodically, the
present paper is very close to the work of L.~Devroye
\cite{Devroye1996} where many examples of mixture representations of
popular probability distributions were discussed from the simulation
point of view. The presented material substantially relies on the
results of \cite{KorolevZeifmanKMJ, Korolevetal2019} and
\cite{LimTeo2009}.

In many situations related to experimental data analysis one often
comes across the following phenomenon: although conventional
reasoning based on the central limit theorem of probability theory
concludes that the expected distribution of observations should be
normal, instead, the statistical procedures expose the noticeable
non-normality of real distributions. Moreover, as a rule, the
observed non-normal distributions are more leptokurtic than the
normal law, having sharper vertices and heavier tails. These
situations are typical in the financial data analysis (see, e. g.,
Chapter 4 in \cite{Shiryaev1998} or Chapter 8 in
\cite{BeningKorolev2002} and the references therein), in
experimental physics (see, e. g., \cite{MeerschaertScheffler2004})
and other fields dealing with statistical analysis of experimental
data. Many attempts were undertaken to explain this
heavy-tailedness. Most significant theoretical breakthrough is
usually associated with the results of B.~Mandelbrot and others
\cite{Mandelbrot1963, Fama1965, Mandelbrot1967} who proposed,
instead of the standard central limit theorem, to use reasoning
based on limit theorems for sums of random summands with infinite
variances (also see \cite{SamorodnitskyTaqqu1994, McCulloch1996})
resulting in non-normal stable laws as heavy-tailed models of the
distributions of experimental data. However, in most cases the key
assumption within this approach, the infiniteness of the variances
of elementary summands, can hardly be believed to hold in practice.
To overcome this contradiction, we consider an extended limit
setting where it may be assumed that the intensity of the flow of
informative events is random resulting in that the number of jumps
up to a certain time in a random-walk-type model or the sample size
is random. We show that in this extended setting, actually,
heavy-tailed scale mixtures of stable laws can also be limit
distributions for sums of a random number of random vectors with
{\it finite} covariance matrices.

The key points of the present paper are:

\begin{itemize}

\item the notion of a scale-mixed multivariate stable distribution is introduced and it is shown that scale-mixed multivariate stable distributions form a special sub-class of multivariate normal scale mixtures;

\item analogs of the muliplication theorem for stable laws are proved for scale-mixed multivariate stable distributions relating these laws with different parameters;

\item some alternative but equivalent definitions are proposed for the generalized multivariate Linnik distributions based on their property to be
scale-mixed multivariate stable distributions;

\item new mixture representations are presented for the generalized multivariate Linnik distributions;

\item a general transfer theorem is proved establishing necessary and sufficient conditions for the convergence of the distributions of sequences of multivariate random vectors with independent random indices (including sums of a random number of random vectors and multivariate statistics constructed from samples with random sizes) to multivariate scale-mixed stable distributions;
\item the property of scale-mixed multivariate stable distributions to be both scale mixtures of a non-trivial multivariate stable distribution and a normal scale mixture is used to obtain necessary and sufficient conditions for the convergence of the distributions of random sums of random vectors to the multivariate generalized Linnik distribution in both cases where the vectors have infinite or finite covariance matrices.

\end{itemize}

This work was partly inspired by the publication of the paper
\cite{JinwenChen2003} in which, based on the results of
\cite{Korolev1994}, a particular case of random sums was considered.
One more reason for writing this work was the recent publication
\cite{SchluterTrede2016}, the authors of which reproduced some
results of \cite{BeningKorolev2005, BeningKorolev2008} and
\cite{Korolev2012} concerning negative binomial sums without citing
these earlier papers.

The paper is organized as follows. In Section 2 we introduce the
main objects of our investigation, the tools we use, and discuss
their properties. Here we present an overview of the properties of
the univariate generalized Mittag-Leffler and generalized Linnik
distributions, introduce the multivariate generalized Linnik
distribution and discuss different approaches to the definition of
the latter. Here we also deduce a simple representation for the
characteristic functions of the scale-mixed multivariate stable
distributions and prove the identifiability of the class of these
mixtures. General properties of scale-mixed multivariate stable
distributions are discussed in Section 3. In this section we also
prove some new mixture representations for the multivariate
generalized Linnik distribution. In Section 4 we, first, prove a
general transfer theorem presenting necessary and sufficient
conditions for the convergence of the distributions of random
sequences with independent random indices (including sums of a
random number of random vectors and multivariate statistics
constructed from samples with random sizes) to scale mixtures of
multivariate elliptically contoured stable distributions. As
particular cases, conditions are obtained for the convergence of the
distributions of scalar normalized random sums of random vectors
with both infinite or finite covariance matrices to scale mixtures
of multivariate stable distributions and their special cases: `pure'
multivariate stable distributions and the multivariate generalized
Linnik distributions. The results of this section extend and refine
those proved in \cite{Korolev1997}.

\section{Preliminaries}

\subsection{Basic notation and definitions}

Let $r\in\mathbb{N}$. We will consider random elements taking values
in the $r$-dimensional Euclidean space $\r^r$. Assume that all the
random variables and random vectors are defined on one and the same
probability space $(\Omega, {\cal A}, {\mbox{\sf P}})$. The
distribution of a random variable $Y$ or an $r$-variate random
vector $\emph{\textbf{Y}}$ with respect to the measure ${\sf P}$
will be denoted ${\cal L}(Y)$ and ${\cal L}(\emph{\textbf{Y}})$,
respectively. The weak convergence, the coincidence of distributions
and the convergence in probability with respect to a specified
probability measure will be denoted by the symbols
$\Longrightarrow$, $\eqd$ and $\pto$, respectively. The product of
{\it independent} random elements will be denoted by the symbol
$\circ$.

A univariate random variable with the standard normal distribution
function $\Phi(x)$ will be denoted $X$,
$$
{\sf
P}(X<x)=\Phi(x)=\frac{1}{\sqrt{2\pi}}\int_{-\infty}^{x}e^{-z^2/2}dz,\
\ \ \ x\in\mathbb{R}.
$$

Let $\Sigma$ be a positive definite $(r\times r)$-matrix. The normal
distribution in $\r^r$ with zero vector of expectations and
covariance matrix $\Sigma$ will be denoted $\mathfrak{N}_{\Sigma}$.
This distribution is defined by its density
$$
\phi(\emph{\textbf{x}})=\frac{\exp\{-\frac12\emph{\textbf{x}}^{\top}\Sigma^{-1}
\emph{\textbf{x}}\}} {(2\pi)^{r/2}|\Sigma|^{1/2}},\ \ \ \
\emph{\textbf{x}}\in\r^r.
$$
The characteristic function
$\mathfrak{f}^{(\emph{\textbf{X}})}(\emph{\textbf{t}})$ of a random
vector $\emph{\textbf{X}}$ such that
$\mathcal{L}(\emph{\textbf{X}})=\mathfrak{N}_{\Sigma}$ has the form
$$
\mathfrak{f}^{(\emph{\textbf{X}})}(\emph{\textbf{t}})\equiv {\sf
E}\exp\{i\emph{\textbf{t}}^{\top}\emph{\textbf{X}}\}=\exp\big\{-{\textstyle\frac12}\emph{\textbf{t}}^{\top}\Sigma\emph{\textbf{t}}\big\},\
\ \ \emph{\textbf{t}}\in\mathbb{R}^r. \eqno(1)
$$

A random variable having the gamma distribution with shape parameter
$r>0$ and scale parameter $\lambda>0$ will be denoted
$G_{r,\lambda}$,
$$
{\sf P}(G_{r,\lambda}<x)=\int_{0}^{x}g(z;r,\lambda)dz,\ \
\text{with}\ \
g(x;r,\lambda)=\frac{\lambda^r}{\Gamma(r)}x^{r-1}e^{-\lambda x},\
x\ge0,
$$
where $\Gamma(r)$ is Euler's gamma-function,
$$
\Gamma(r)=\int_{0}^{\infty}x^{r-1}e^{-x}dx, \ \ \ \ r>0.
$$
In this notation, obviously, $G_{1,1}$ is a random variable with the
standard exponential distribution: ${\sf
P}(G_{1,1}<x)=\big[1-e^{-x}\big]{\bf 1}(x\ge0)$ (here and in what
follows ${\bf 1}(A)$ is the indicator function of a set $A$).

Let $D_{r}$ be a random variable with the {\it one-sided exponential
power distribution} defined by the density
$$
f^{(D)}_{r}(x)=\frac{1}{r\Gamma(1/r)}e^{-x^{1/r}},\ \ \ x\ge0.
$$
It is easy to make sure that $G_{r,1}\eqd D_{r}^{1/r}$.

The gamma distribution is a particular representative of the class
of generalized gamma distributions (GG distributions), that was
first described in \cite{Stacy1962} as a special family of lifetime
distributions containing both gamma and Weibull distributions. A
{\it generalized gamma $($GG$)$ distribution} is the absolutely
continuous distribution defined by the density
$$
\overline{g}(x;r,\alpha,\lambda)=\frac{|\alpha|\lambda^r}{\Gamma(r)}x^{\alpha
r-1}e^{-\lambda x^{\alpha}},\ \ \ \ x\ge0,
$$
with $\alpha\in\mathbb{R}$, $\lambda>0$, $r>0$. A random variable
with the density $\overline{g}(x;r,\alpha,\lambda)$ will be denoted
$\overline{G}_{r,\alpha,\lambda}$. It is easy to see that
$$
\overline{G}_{r,\alpha,\mu}\eqd G_{r,\mu}^{1/\alpha}\eqd
\mu^{-1/\alpha}G_{r,1}^{1/\alpha}\eqd\mu^{-1/\alpha}\overline{G}_{r,\alpha,1}.
$$

Let $\gamma>0$. The distribution of the random variable
$W_{\gamma}$:
$$
{\sf
P}\big(W_{\gamma}<x\big)=\big[1-e^{-x^{\gamma}}\big]\mathbf{1}(x\ge
0),
$$
is called the {\it Weibull distribution} with shape parameter
$\gamma$. It is obvious that $W_1$ is the random variable with the
standard exponential distribution: ${\sf
P}(W_1<x)=\big[1-e^{-x}\big]{\bf 1}(x\ge0)$. The Weibull
distribution is a particular case of GG distributions corresponding
to the density $\overline{g}(x;1,\gamma,1)$. It is easy to see that
$W_1^{1/\gamma}\eqd W_{\gamma}$. Moreover, if $\gamma>0$ and
$\gamma'>0$, then ${\sf P}(W_{\gamma'}^{1/\gamma}\ge x)={\sf
P}(W_{\gamma'}\ge x^{\gamma})=e^{-x^{\gamma\gamma'}}={\sf
P}(W_{\gamma\gamma'}\ge x)$, $x\ge 0$, that is, for any $\gamma>0$
and $\gamma'>0$
$$
W_{\gamma\gamma'}\eqd W_{\gamma'}^{1/\gamma}.\eqno(2) % (3)
$$

\smallskip

In the paper \cite{Gleser1989} it was shown that any gamma
distribution with shape parameter no greater than one is mixed
exponential. Namely, the density $g(x;r,\mu)$ of a gamma
distribution with $0<r<1$ can be represented as
$$
g(x;r,\mu)=\int_{0}^{\infty}ze^{-zx}p(z;r,\mu)dz,
$$
where
$$
p(z;r,\mu)=\frac{\mu^r}{\Gamma(1-r)\Gamma(r)}\cdot\frac{\mathbf{1}(z\ge\mu)}{(z-\mu)^rz}.\eqno(3)
$$
Moreover, a gamma distribution with shape parameter $r>1$ cannot be
represented as a mixed exponential distribution.

In \cite{Korolev2017} it was proved that if $r\in(0,1)$, $\mu>0$ and
$G_{r,\,1}$ and $G_{1-r,\,1}$ are independent gamma-distributed
random variables, then the density $p(z;r,\mu)$ defined by $(3)$
corresponds to the random variable
$$
Z_{r,\mu}=\frac{\mu(G_{r,\,1}+G_{1-r,\,1})}{G_{r,\,1}}\eqd\mu
Z_{r,1}\eqd\mu\big(1+\textstyle{\frac{1-r}{r}}V_{1-r,r}\big),
$$
where $V_{1-r,r}$ is the random variable with the Snedecor--Fisher
distribution defined by the probability density
$$
q(x;1-r,r)=\frac{(1-r)^{1-r}r^r}{\Gamma(1-r)\Gamma(r)}
\cdot\frac{1}{x^{r}[r+(1-r)x]},\ \ \ x\ge0.
$$
In other words, if $r\in(0,1)$, then
$$
G_{r,\,\mu}\eqd W_1\circ Z_{r,\,\mu}^{-1}.\eqno(4)
$$

\subsection{Stable distributions}

Any random variable that has the univariate strictly stable
distribution with the characteristic exponent $\alpha$ and shape
parameter $\theta$ corresponding to the characteristic function
$$
\mathfrak{f}_{\alpha,\theta}(t)=\exp\big\{-|t|^{\alpha}\exp\{-{\textstyle\frac12}i\pi\theta\alpha\,\mathrm{sign}t\}\big\},\
\ \ \ t\in\r,\eqno(5) % (2)
$$
with $0<\alpha\le2$, $|\theta|\le\min\{1,\frac{2}{\alpha}-1\}$, will
be denoted $S(\alpha,\theta)$ (see, e. g., \cite{Zolotarev1986}).
For definiteness, $S(1,\,1)=1$.

From (5) it follows that the characteristic function of a symmetric
($\theta=0$) strictly stable distribution has the form
$$
\mathfrak{f}_{\alpha,0}(t)=e^{-|t|^{\alpha}},\ \ \ t\in\r. \eqno(6) % (3)
$$
From (6) it is easy to see that $S(2,0)\eqd\sqrt{2}X$.

Let $r\in\mathbb{N}$, $\alpha\in(0,2]$. An $r$-variate random vector
$\emph{\textbf{S}}(\alpha,\Sigma,0)$ is said to have the (centered)
elliptically contoured stable distribution
$\mathfrak{S}_{\alpha,\Sigma,0}$ with characteristic exponent
$\alpha$, if its characteristic function
$\mathfrak{f}_{\alpha,\Sigma,0}(\emph{\textbf{t}})$ has the form
$$
\mathfrak{f}_{\alpha,\Sigma,0}(\emph{\textbf{t}})\equiv{\sf
E}\exp\{i\emph{\textbf{t}}^{\top}\emph{\textbf{S}}(\alpha,\Sigma,0)\}=
\exp\{-(\emph{\textbf{t}}^{\top}\Sigma\emph{\textbf{t}})^{\alpha/2}\},\
\ \ \emph{\textbf{t}}\in\r^r.
$$
It is easy to see that
$\mathfrak{S}_{2,\Sigma,0}=\mathfrak{N}_{2\Sigma}$.

Univariate stable distributions are popular examples of heavy-tailed
distributions. Their moments of orders $\delta\ge\alpha$ do not
exist (the only exception is the normal law corresponding to
$\alpha=2$), and if $0<\delta<\alpha$, then
$$
{\sf E}|S(\alpha,0)|^{\delta}=\frac{2^{\delta}}{\sqrt{\pi}}\cdot
\frac{\Gamma(\frac{\delta+1}{2})\Gamma(1-\frac{\delta}{\alpha})}{\Gamma(\frac{2}{\delta}-1)}\eqno(7)
$$
(see, e. g., \cite{Korolev2016}). Stable laws and only they can be
limit distributions for sums of a non-random number of independent
identically distributed random variables with infinite variance
under linear normalization.

Let $0<\alpha\le 1$. By $S(\alpha,1)$ we will denote a positive
random variable with the one-sided stable distribution corresponding
to the characteristic function
$$
\mathfrak{f}_{\alpha}(t)=\exp\big\{-|t|^{\alpha}\exp\{-{\textstyle\frac12}i\pi\alpha\,\mathrm{sign}t\}\big\},\
\ \ \ t\in\r.
$$
The Laplace--Stieltjes transform $\psi^{(S)}_{\alpha,1}(s)$ of the
random variable $S(\alpha,1)$ has the form
$$
\psi^{(S)}_{\alpha,1}(s)\equiv{\sf E}\exp\{-sS(\alpha,1)\}
=e^{-s^{\alpha}},\ \ \ s>0.
$$
The moments of orders $\delta\ge\alpha$ of the random variable
$S(\alpha,1)$ are infinite and for $0<\delta<\alpha$ we have
$$
{\sf
E}S^{\delta}(\alpha,1)=\frac{2^{\delta}\Gamma(1-\frac{\delta}{\alpha})}{\Gamma(1-\delta)}\eqno(8)
$$
(see, e. g., \cite{Korolev2016}). For more details see
\cite{Zolotarev1986} or \cite{SamorodnitskyTaqqu1994}.

It is known that if $0<\alpha\le1$ and $0<\alpha'\le1$, then
$$
S(\alpha\alpha',1)\eqd S^{1/\alpha}(\alpha',1)\circ
S(\alpha,1),\eqno(9)  %(19)
$$
see Corollary 1 to Theorem 3.3.1 in \cite{Zolotarev1986}.

Let $\alpha\in(0,2]$. It is known that, if $\emph{\textbf{X}}$ is a
random vector such that
$\mathcal{L}(\emph{\textbf{X}})=\mathfrak{N}_{\Sigma}$ independent
of the random variable $S(\alpha/2,1)$, then
$$
\emph{\textbf{S}}(\alpha,\Sigma,0)\eqd S^{1/2}(\alpha/2,1)\circ
\emph{\textbf{S}}(2,\Sigma,0)\eqd \sqrt{2S(\alpha/2,1)}\circ
\emph{\textbf{X}}\eqno(10) %(20)
$$
(see Proposition 2.5.2 in \cite{SamorodnitskyTaqqu1994}).

It is easy to make sure that if $0<\alpha\le 2$ and $0<\alpha'\le1$,
then
$$
\emph{\textbf{S}}(\alpha\alpha',\Sigma,0)\eqd
S^{1/\alpha}(\alpha',1)\circ\emph{\textbf{S}}(\alpha,\Sigma,0).\eqno(11)%(23)
$$
Indeed, from (9) and (10) we have
$$
\emph{\textbf{S}}(\alpha\alpha',\Sigma,0)\eqd
\sqrt{S(\alpha\alpha'/2,1)}\circ \emph{\textbf{S}}(2,\Sigma,0)\eqd
\sqrt{2S(\alpha\alpha'/2,1)}\circ \emph{\textbf{X}}\eqd
$$
$$
\eqd\sqrt{S^{2/\alpha}(\alpha',1)}\circ\sqrt{2S(\alpha/2,1)}\circ\emph{\textbf{X}}\eqd
S^{1/\alpha}(\alpha',1)\circ\emph{\textbf{S}}(\alpha,\Sigma,0).
$$

If $\alpha=2$, then (11) turns into (10).

\subsection{Scale mixtures of multivariate distributions}

Let $U$ be a nonnegative random variable. The symbol ${\sf
E}\mathfrak{N}_{U\Sigma}(\cdot)$ will denote the distribution which
for each Borel set $A$ in $\r^r$ is defined as
$$
{\sf
E}\mathfrak{N}_{U\Sigma}(A)=\int_{0}^{\infty}\mathfrak{N}_{u\Sigma}(A)d{\sf
P}(U<u).
$$

It is easy to see that if $\emph{\textbf{X}}$ is a random vector
such that $\mathcal{L}(\emph{\textbf{X}})=\mathfrak{N}_{\Sigma}$,
then ${\sf
E}\mathfrak{N}_{U\Sigma}=\mathcal{L}(\sqrt{U}\circ\emph{\textbf{X}})$.

In this notation, relation (10) can be written as
$$
\mathfrak{S}_{\alpha,\Sigma,0}={\sf
E}\mathfrak{N}_{2S(\alpha/2,1)\Sigma}.\eqno(12)   %(22)
$$

By analogy, the symbol ${\sf
E}\mathfrak{S}_{\alpha,U^{2/\alpha}\Sigma,0}$ will denote the
distribution that for each Borel set $A$ in $\r^r$ is defined as
$$
{\sf
E}\mathfrak{S}_{\alpha,U^{2/\alpha}\Sigma,0}(A)=\int_{0}^{\infty}\mathfrak{S}_{\alpha,u^{2/\alpha}\Sigma,0}(A)d{\sf
P}(U<u).
$$
The characteristic function corresponding to the distribution ${\sf
E}\mathfrak{S}_{\alpha,0,U^{2/\alpha}\Sigma}$ has the form
$$
\int_{0}^{\infty}\exp\big\{-\big(\emph{\textbf{t}}^{\top}(u^{2/\alpha}\Sigma)\emph{\textbf{t}}\big)^{\alpha/2}\big\}
d{\sf
P}(U<u)=\int_{0}^{\infty}\exp\big\{-\big((u^{1/\alpha}\emph{\textbf{t}})^{\top}\Sigma(u^{1/\alpha}\emph{\textbf{t}})\big)^{\alpha/2}\big\}
d{\sf P}(U<u)=
$$
$$
={\sf
E}\exp\big\{i\emph{\textbf{t}}^{\top}U^{1/\alpha}\circ\emph{\textbf{S}}(\alpha,\Sigma,0)\big\},\
\ \ \ \emph{\textbf{t}}\in\r^r,\eqno(13)  %(24)
$$
where the random variable $U$ is independent of the random vector
$\emph{\textbf{S}}(\alpha,\Sigma,0)$, that is, the distribution
${\sf E}\mathfrak{S}_{\alpha,U^{2/\alpha}\Sigma,0}$ corresponds to
the product
$U^{1/\alpha}\circ\cdot\emph{\textbf{S}}(\alpha,\Sigma,0)$.

Let ${\cal U}$ be the set of all nonnegative random variables. Now
consider an auxiliary statement dealing with the identifiability of
the family of distributions $\{{\sf
E}\mathfrak{S}_{\alpha,U^{2/\alpha}\Sigma,0}:\ U\in{\cal U}\}$.

\smallskip

{\sc Lemma 1}. {\it Whatever a nonsingular positive definite matrix
$\Sigma$ is, the family $\{{\sf
E}\mathfrak{S}_{\alpha,U^{2/\alpha}\Sigma,0}:\ U\in{\cal U}\}$ is
identifiable in the sense that if $U_1\in{\cal U}$, $U_2\in{\cal U}$
and
$$
{\sf E}\mathfrak{S}_{\alpha,U_1^{2/\alpha}\Sigma,0}(A)={\sf
E}\mathfrak{S}_{\alpha,U_2^{2/\alpha}\Sigma,0}(A)\eqno(14)   %(4)
$$
for any set $A\in{\cal B}(\r^r)$, then $U_1\eqd U_2$.}

\smallskip

The $\,$ {\sc proof} $\,$ of this lemma is very simple. If
$U\in{\cal U}$, then it follows from (13) that the characteristic
function $\mathfrak{v}^{(U)}_{\alpha,\Sigma}(\emph{\textbf{t}})$
corresponding to the distribution ${\sf
E}\mathfrak{S}_{\alpha,U^{2/\alpha}\Sigma,0}$ has the form
$$
\mathfrak{v}^{(U)}_{\alpha,\Sigma}(\emph{\textbf{t}})=
\int_{0}^{\infty}\exp\big\{-\big(\emph{\textbf{t}}^{\top}(u^{2/\alpha}\Sigma)\emph{\textbf{t}}\big)^{\alpha/2}\big\}
d{\sf P}(U<u)=
$$
$$
= \int_{0}^{\infty}\exp\{-us\}d{\sf P}(U<u),\ \
s=(\emph{\textbf{t}}^{\top}\Sigma\emph{\textbf{t}})^{\alpha/2},\ \ \
 \ \emph{\textbf{t}}\in\r^r, \eqno(15)    %(5)
$$
But on the right-hand side of (15) there is the Laplace--Stieltjes
transform of the random variable $U$. From (14) it follows that
$\mathfrak{v}^{(U_1)}_{\alpha,\Sigma}(\emph{\textbf{t}})\equiv
\mathfrak{v}^{(U_2)}_{\alpha,\Sigma}(\emph{\textbf{t}})$ whence by
virtue of (15) the Laplace--Stieltjes transforms of the random
variables $U_1$ and $U_2$ coincide, whence, in turn, it follows that
$U_1\eqd U_2$. The lemma is proved.

\smallskip

{\sc Remark 1.} When proving Lemma 1 we established a simple but
useful by-product result: if $\psi^{(U)}(s)$ is the
Laplace--Stieltjes transform of the random variable $U$, then the
characteristic function
$\mathfrak{v}^{(U)}_{\alpha,\Sigma}(\emph{\textbf{t}})$
corresponding to the distribution ${\sf
E}\mathfrak{S}_{\alpha,U^{2/\alpha}\Sigma,0}$ has the form
$$
\mathfrak{v}^{(U)}_{\alpha,\Sigma}(\emph{\textbf{t}})=
\psi^{(U)}\big((\emph{\textbf{t}}^{\top}\Sigma\emph{\textbf{t}})^{\alpha/2}\big),\
\ \emph{\textbf{t}}\in\r^r. \eqno(16)  %(6)
$$

\smallskip

Let $\emph{\textbf{X}}$ be a random vector such that
$\mathcal{L}(\emph{\textbf{X}})=\mathfrak{N}_{\Sigma}$ with some
positive definite $(r\times r)$-matrix $\Sigma$. Define the
multivariate Laplace distribution as
$\mathcal{L}\big(\sqrt{2W_1}\circ \emph{\textbf{X}}\big)={\sf
E}\mathfrak{N}_{2W_1\Sigma}$. The random vector with this
multivariate Laplace distribution will be denoted
$\Lambda_{\Sigma}$. It is well known that the Laplace---Stieltjes
transform $\psi^{(W_1)}(s)$ of the random variable $W_1$ with the
exponential distribution has the form
$$
\psi^{(W_1)}(s)=(1+s)^{-1},\ \ \ \ s>0.\eqno(17)
$$
Hence, in accordance with (17) and Remark 1, the characteristic
function $\mathfrak{f}^{(\Lambda)}_{\Sigma}(\emph{\textbf{t}})$ of
the random variable $\Lambda_{\Sigma}$ has the form
$$
\mathfrak{f}^{(\Lambda)}_{\Sigma}(\emph{\textbf{t}})=
\psi^{(W_1)}\big(\emph{\textbf{t}}^{\top}\Sigma\emph{\textbf{t}}\big)=
\big(1+\emph{\textbf{t}}^{\top}\Sigma\emph{\textbf{t}}\big)^{-1},\ \
\ \emph{\textbf{t}}\in\mathbb{R}^r.
$$

\subsection{The generalized Mittag-Leffler distribution}

The probability distribution of a nonnegative random variable
$M_{\delta}$ whose Laplace transform is
$$
\psi^{(M)}_{\delta}(s)\equiv {\sf E}e^{-sM_{\delta}}=\big(1+\lambda
s^{\delta}\big)^{-1},\ \ \
s\ge0,\eqno(18)   %(4)
$$
where $\lambda>0$, $0<\delta\le1$, is called {\it the Mittag-Leffler
distribution}. For simplicity, in what follows we will consider the
standard scale case and assume that $\lambda=1$.

The origin of the term {\it Mittag-Leffler distribution} is due to
that the probability density corresponding to Laplace transform (18)
has the form
$$
f_{\delta}^{(M)}(x)=\frac{1}{x^{1-\delta}}\sum\nolimits_{n=0}^{\infty}\frac{(-1)^nx^{\delta
n}}{\Gamma(\delta n+1)}=-\frac{d}{dx}E_{\delta}(-x^{\delta}),\ \ \
x\ge0,\eqno(19)     %(5)
$$
where $E_{\delta}(z)$ is the Mittag-Leffler function with index
$\delta$ that is defined as the power series
$$
E_{\delta}(z)=\sum\nolimits_{n=0}^{\infty}\frac{z^n}{\Gamma(\delta
n+1)},\ \ \ \delta>0,\ z\in\mathbb{Z}.
$$
%The distribution function corresponding to density (19) will be
%denoted $F_{\delta}^{M}(x)$.

With $\delta=1$, the Mittag-Leffler distribution turns into the
standard exponential distribution, that is, $F_1^{M}(x)=
[1-e^{-x}]\mathbf{1}(x\ge 0)$, $x\in\mathbb{R}$. But with $\delta<1$
the Mittag-Leffler distribution density has the heavy power-type
tail: from the well-known asymptotic properties of the
Mittag-Leffler function it can be deduced that if $0<\delta<1$, then
$$
f_\delta^{(M)}(x)\sim \frac{\sin(\delta\pi)\Gamma(\delta+1)}{\pi
x^{\delta+1}}
$$
as $x\to\infty$, see, e. g., \cite{Kilbas2014}.

It is well-known that the Mittag-Leffler distribution is
geometrically stable. This means that if $X_1,X_2,\ldots$ are
independent random variables whose distributions belong to the
domain of attraction of a one-sided $\alpha$-strictly stable law
$\mathcal{L}(S(\alpha,1))$ and $\nb_{1,\,p}$ is the random variable
independent of $X_1,X_2,\ldots$ and having the geometric
distribution
$$
{\sf P}(\nb_{1,\,p}=n)=p(1-p)^{n-1},\ \ \ n=1,2,\ldots,\ \ \
p\in(0,1),\eqno(20)
$$
then for each $p\in(0,1)$ there exists a constant $a_p>0$ such that
$a_p\big(X_1+\ldots+X_{\nb_{1,\,p}}\big)\Longrightarrow M_{\delta}$
as $p\to 0$, see, e. g., \cite{KlebanovRachev1996}. The history of
the Mittag-Leffler distribution is discussed in
\cite{KorolevZeifmanKMJ}. For more details see e. g.,
\cite{KorolevZeifman2016, KorolevZeifmanKMJ} and the references
therein. The Mittag-Leffler distributions are of serious theoretical
interest in the problems related to thinned (or rarefied)
homogeneous flows of events such as renewal processes or anomalous
diffusion or relaxation phenomena, see \cite{WeronKotulski1996,
GorenfloMainardi2006} and the references therein.

Let $\nu>0$, $\delta\in(0,1]$. The distribution of a nonnegative
random variable $M_{\delta,\,\nu}$ defined by the Laplace--Stieltjes
transform
$$
\psi^{(M)}_{\delta,\,\nu}(s)\equiv {\sf
E}e^{-sM_{\delta,\,\nu}}=\big(1+s^{\delta}\big)^{-\nu},\ \ \
s\ge0,\eqno(21)
$$
is called the {\it generalized Mittag-Leffler distribution}, see
\cite{MathaiHaubold2011, Joseetal} and the references therein.
Sometimes this distribution is called the {\it Pillai distribution}
\cite{Devroye1996}, although in the original paper \cite{Pillai1985}
R. Pillai called it {\it semi-Laplace}. In the present paper we will
keep to the first term {\it generalized Mittag-Leffler
distribution}.

The properties of univariate generalized Mittag-Leffler distribution
are discussed in \cite{Mathai2010, MathaiHaubold2011, Joseetal,
Korolevetal2018}. In particular, it is well known that if
$\delta\in(0,1]$ and $\nu>0$, then
$$
M_{\delta,\,\nu}\eqd
S(\delta,\,1)\circ\overline{G}_{\nu,\,\delta,\,1}\eqd
S(\delta,\,1)\circ G^{1/\delta}_{\nu,1}\eqno(22)
$$
(see \cite{MathaiHaubold2011, Joseetal}). If $\beta\ge\delta$, then
the moments of order $\beta$ of the random variable
$M_{\delta,\,\nu}$ are infinite, and if $0<\beta<\delta<1$, then
$$
{\sf
E}M_{\delta,\,\nu}^{\beta}=\frac{\Gamma(1-\frac{\beta}{\delta})\Gamma(\nu+\frac{\beta}{\delta})}{\Gamma(1-\beta)\Gamma(\nu)},
$$
see \cite{Korolevetal2018}.

In \cite{Korolevetal2018} it was demonstrated that the generalized
Mittag-Leffler distribution can be represented as a scale mixture of
`ordinary' Mittag-Leffler distributions: if $\nu\in(0,1]$ and
$\delta\in(0,1]$, then
$$
M_{\delta,\,\nu}\eqd Z_{\nu,1}^{-1/\delta}\circ M_{\delta}.\eqno(23)
$$
In \cite{Korolevetal2018} it was also shown that any generalized
Mittag-Leffler distribution is a scale mixture a one-sided stable
law with any greater characteristic parameter, the mixing
distribution being the generalized Mittag-Leffler law: if
$\delta\in(0,1]$, $\delta'\in(0,1)$ and $\nu>0$, then
$$
M_{\delta\delta',\,\nu}\eqd S(\delta,1)\circ
M_{\delta',\nu}^{1/\delta}.\eqno(24)
$$

\subsection{The generalized Linnik distributions}

In 1953 Yu. V. Linnik \cite{Linnik1953} introduced a class of
symmetric distributions whose characteristic functions have the form
$$
\mathfrak{f}^{(L)}_{\alpha}(t)=\big(1+|t|^{\alpha}\big)^{-1},\ \ \
t\in\mathbb{R},\eqno(25)   %(1)
$$
where $\alpha\in(0,2]$. The distributions with the characteristic
function (25) are traditionally called the {\it Linnik
distributions}. Although sometimes the term {\it $\alpha$-Laplace
distributions} \cite{Pillai1985} is used, we will use the first term
which has already become conventional. If $\alpha=2$, then the
Linnik distribution turns into the Laplace distribution
corresponding to the density
$$
f^{(\Lambda)}(x)=\textstyle{\frac12}e^{-|x|},\ \ \
x\in\mathbb{R}.\eqno(26)
$$
A random variable with density (26) will be denoted $\Lambda$. A
random variable with the Linnik distribution with parameter $\alpha$
will be denoted $L_{\alpha}$.

Perhaps, most often Linnik distributions are recalled as examples of
symmetric geometric stable distributions. This means that if
$X_1,X_2,\ldots$ are independent random variables whose
distributions belong to the domain of attraction of an
$\alpha$-strictly stable symmetric law and $\nb_{1,\,p}$ is the
random variable independent of $X_1,X_2,\ldots$ and having the
geometric distribution (20), then for each $p\in(0,1)$ there exists
a constant $a_p>0$ such that
$a_p\big(X_1+\ldots+X_{\nb_{1,\,p}}\big)\Longrightarrow L_{\alpha}$
as $p\to 0$, see, e. g., \cite{Bunge1996} or
\cite{KlebanovRachev1996}.

The properties of the Linnik distributions were studied in many
papers. We should mention \cite{Laha1961, Devroye1990,
KotzOstrovskiiHayfavi1995a, KotzOstrovskiiHayfavi1995b} and other
papers, see the survey in \cite{KorolevZeifmanKMJ}.

In \cite{Devroye1990} and \cite{KorolevZeifmanKMJ} it was
demonstrated that
$$
L_{\alpha}\eqd W_1^{1/\alpha}\circ S(\alpha,0)\eqd \sqrt{2M_{\alpha/2}}\circ X,\eqno(27)  %(25)
$$
where the random variable $M_{\alpha/2}$ has the Mittag-Leffler
distribution with parameter $\alpha/2$.

The multivariate Linnik distribution was introduced by D. N.
Anderson in \cite{Anderson1992} where it was proved that the
function
$$
\mathfrak{f}^{(L)}_{\alpha,\Sigma}(\emph{\textbf{t}})=\big[1+
(\emph{\textbf{t}}^{\top}\Sigma\emph{\textbf{t}})^{\alpha/2}\big]^{-1},\
\ \
\emph{\textbf{t}}\in\r^r,\ \ \alpha\in(0,2),\eqno(28)   %(23)
$$
is the characteristic function of an $r$-variate probability
distribution, where $\Sigma$ is a positive definite $(r\times
r)$-matrix. In \cite{Anderson1992} the distribution corresponding to
the characteristic function (28) was called {\it the $r$-variate
Linnik distribution}. For the properties of the multivariate Linnik
distributions see \cite{Anderson1992, Ostrovskii1995}.

The $r$-variate Linnik distribution can also be defined in another
way. Let $\emph{\textbf{X}}$ be a random vector such that
$\mathcal{L}(\emph{\textbf{X}})=\mathfrak{N}_{\Sigma}$, where
$\Sigma$ is a positive definite $(r\times r)$-matrix, independent of
the random variable $M_{\alpha/2}$. By analogy with (27) introduce
the random vector $\emph{\textbf{L}}_{\alpha,\Sigma}$ as
$$
\emph{\textbf{L}}_{\alpha,\Sigma}=\sqrt{2M_{\alpha/2}}\circ\emph{\textbf{X}}.
$$
Then, in accordance with what has been said in Section 2.3,
$$
\mathcal{L}(\emph{\textbf{L}}_{\alpha,\Sigma})={\sf
E}\mathfrak{N}_{2M_{\alpha/2}\Sigma}.\eqno(29)   %(26)
$$
The distribution (29) will be called {\it the $($centered$)$
elliptically contoured multivariate Linnik distribution}.

Using Remark 1 we can easily make sure that the two definitions of
the multivariate Linnik distribution coincide. Indeed, with the
account of (18), according to Remark 1, the characteristic function
of the random vector $\emph{\textbf{L}}_{\alpha,\Sigma}$ defined by
(29) has the form
$$
{\sf
E}\exp\{i\emph{\textbf{t}}^{\top}\emph{\textbf{L}}_{\alpha,\Sigma}\}=
\psi^{(M)}_{\alpha/2}\big(\emph{\textbf{t}}^{\top}\Sigma\emph{\textbf{t}}\big)=
\big[1+(\emph{\textbf{t}}^{\top}\Sigma\emph{\textbf{t}})^{\alpha/2}\big]^{-1}=\mathfrak{f}^{(L)}_{\alpha,\Sigma}(\emph{\textbf{t}}),\
\ \emph{\textbf{t}}\in\r^r,
$$
that coincides with Anderson's definition (28).

Based on (27), one more equivalent definition of the multivariate
Linnik distribution can be proposed. Namely, let
$\emph{\textbf{L}}_{\alpha,\Sigma}$ be an $r$-variate random vector
such that
$$
\emph{\textbf{L}}_{\alpha,\Sigma}=W_1^{1/\alpha}\circ\emph{\textbf{S}}(\alpha,\Sigma,0).\eqno(30)
$$
In accordance with (17) and Remark 1 the characteristic function of
the random vector $\emph{\textbf{L}}_{\alpha,\Sigma}$ defined by
(30) again has the form
$$
{\sf
E}\exp\{i\emph{\textbf{t}}^{\top}\emph{\textbf{L}}_{\alpha,\Sigma}\}=
\psi^{(W_1)}\big((\emph{\textbf{t}}^{\top}\Sigma\emph{\textbf{t}})^{\alpha/2}\big)=
\big[1+(\emph{\textbf{t}}^{\top}\Sigma\emph{\textbf{t}})^{\alpha/2}\big]^{-1}=\mathfrak{f}^{(L)}_{\alpha,\Sigma}(\emph{\textbf{t}}),\
\ \emph{\textbf{t}}\in\r^r.
$$

The definitions (29) and (30) open the way to formulate limit
theorems stating that the multivariate Linnik distribution can not
only be limiting for geometric random sums of independent
identically distributed random vectors with infinite second moments
\cite{KozubowskiRachev1999}, but it also can be limiting for random
sums of independent random vectors with finite covariance matrices.

In \cite{Pakes1998}, Pakes showed that the probability distributions
known as {\it generalized Linnik distributions} which have
characteristic functions
$$
\mathfrak{f}^{(L)}_{\alpha,\nu}(t)=\big(1+|t|^{\alpha}\big)^{-\nu},\
\ \ t\in\mathbb{R},\ 0<\alpha\le 2,\ \nu>0,\eqno(31)    %(3)
$$
play an important role in some characterization problems of
mathematical statistics. The class of probability distributions
corresponding to characteristic function (31) have found some
interesting properties and applications, see \cite{Anderson1993,
BaringhausGrubel1997, Devroye1990, Jayakumar1995,
KotzOstrovskii1996, Kotz2001, Kozubowski1998, Lin1998} and related
papers. In particular, they are good candidates to model financial
data which exhibits high kurtosis and heavy tails
\cite{MittnikRachev1993}.

%Here we concentrate our attention at the symmetric case $\theta=0$.
Any random variable with the characteristic function (31) will be
denoted $L_{\alpha,\nu}$.

Recall some results containing mixture representations for the
generalized Linnik distribution. The following well-known result is
due to Devroye \cite{Devroye1990} and Pakes \cite{Pakes1998} who
showed that
$$
L_{\alpha,\nu}\eqd S_{\alpha,0}\circ G_{\nu,1}^{1/\alpha}\eqd
S_{\alpha,0}\cdot\overline{G}_{\nu,\alpha,1}\eqno(32)
$$
for any $\alpha\in(0,2]$ and $\nu>0$.

It is well known that
$$
{\sf E}G_{\nu,1}^{\gamma}=\frac{\Gamma(\nu+\gamma)}{\Gamma(\nu)}
$$
for $\gamma>-\nu$. Hence, for $0\le\beta<\alpha$ from (7) and (32)
we obtain
$$
{\sf E}|L_{\alpha,\nu}|^{\beta}={\sf
E}|S_{\alpha,0}|^{\beta}\cdot{\sf
E}G_{\nu,1}^{\beta/\alpha}=\frac{2^{\beta}}{\sqrt{\pi}}\cdot
\frac{\Gamma(\frac{\beta+1}{2})\Gamma(1-\frac{\beta}{\alpha})\Gamma(\nu+\frac{\beta}{\alpha})}{\Gamma(\frac{2}{\beta}-1)\Gamma(\nu)}.
$$

Generalizing and improving some results of \cite{Pakes1998} and
\cite{LimTeo2009}, with the account of (22) in
\cite{Korolevetal2019} it was demonstrated that for $\nu>0$ and
$\alpha\in(0,2]$
$$
L_{\alpha,\nu}\eqd X\circ\sqrt{2S(\alpha/2,1)}\circ
G_{\nu,1}^{1/\alpha}\eqd
X\circ\sqrt{2S(\alpha/2,1)\circ\overline{G}_{\nu,\alpha/2,1}}\eqd
X\circ\sqrt{2M_{\alpha/2,\,\nu}}.\eqno(33)
   %(16)
$$
that is, the generalized Linnik distribution is a normal scale
mixture with the generalized Mittag-Leffler mixing distribution.

It is easy to see that for any $\alpha>0$ and $\alpha'>0$
$$
\overline{G}_{\nu,\alpha\alpha',1}\eqd
G_{\nu,1}^{1/\alpha\alpha'}\eqd
(G_{\nu,1}^{1/\alpha'})^{1/\alpha}\eqd\overline{G}_{\nu,\alpha',1}^{1/\alpha}.\eqno(34)
$$
Therefore, for $\alpha\in(0,2]$, $\alpha'\in(0,1)$ and $\nu>0$ using
(32) and the univariate version of (11) we obtain the following
chain of relations:
$$
L_{\alpha\alpha',\,\nu}\eqd S(\alpha\alpha',0)\circ
G_{\nu,1}^{1/\alpha\alpha'}\eqd S(\alpha,0)\circ
S^{1/\alpha}(\alpha',1)\circ G_{\nu,1}^{1/\alpha\alpha'}\eqd
$$
$$
\eqd S(\alpha,0)\circ
\big(S(\alpha',1)\overline{G}_{\nu,\alpha',1}\big)^{1/\alpha}\eqd
S(\alpha,0)\circ M_{\alpha',\nu}^{1/\alpha}.
$$
Hence, the following statement, more general than (33), holds
representing the generalized Linnik distribution as a scale mixture
of a symmetric stable law with any greater characteristic parameter,
the mixing distribution being the generalized Mittag-Leffler law: if
$\alpha\in(0,2]$, $\alpha'\in(0,1)$ and $\nu>0$, then
$$
L_{\alpha\alpha',\,\nu}\eqd S(\alpha,0)\circ
M_{\alpha',\nu}^{1/\alpha}.\eqno(35)
$$

Now let $\nu\in(0,1]$. From (32) and (4) it follows that
$$
L_{\alpha,\nu}\eqd S(\alpha,0)\circ G_{\nu,1}^{1/\alpha}\eqd
S(\alpha,0)\circ W_1^{1/\alpha}\circ Z_{\nu,1}^{-1/\alpha}\eqd
L_{\alpha}\circ Z_{\nu,1}^{-1/\alpha}
$$
yielding the following relation proved in \cite{Korolevetal2019}: if
$\nu\in(0,1]$ and $\alpha\in(0,2]$, then
$$
L_{\alpha,\nu}\eqd L_{\alpha}\cdot Z_{\nu,1}^{-1/\alpha}.\eqno(36)
$$
In other words, with $\nu\in(0,1]$ and $\alpha\in(0,2]$, the
generalized Linnik distribution is a scale mixture of `ordinary'
Linnik distributions. In the same paper the representation of the
generalized Linnik distribution via the Laplace and `ordinary'
Mittag-Leffler distributions was obtained.

For $\delta\in(0,1]$ denote
$$
R_{\delta}=\frac{S(\delta,1)}{S'(\delta,1)},
$$
where $S(\delta,1)$ and $S'(\delta,1)$ are independent random
variables with one and the same one-sided stable distribution with
the characteristic exponent $\delta$. In \cite{KorolevZeifmanKMJ} it
was shown that the probability density $f^{(R)}_{\delta}(x)$ of the
ratio $R_{\delta}$ of two independent random variables with one and
the same one-sided strictly stable distribution with parameter
$\delta$ has the form
$$
f^{(R)}_{\delta}(x)=
\frac{\sin(\pi\delta)x^{\delta-1}}{\pi[1+x^{2\delta}+2x^{\delta}\cos(\pi\delta)]},\
\ \ x>0.
$$
In \cite{Korolevetal2019} it was proved that if $\nu\in(0,1]$ and
$\alpha\in(0,2]$, then
$$
L_{\alpha,\nu}\eqd X\circ
Z_{\nu,1}^{-1/\alpha}\circ\sqrt{2M_{\alpha/2}}\eqd \Lambda\circ
Z_{\nu,1}^{-1/\alpha}\circ\sqrt{R_{\alpha/2}}.\eqno(37)
$$
So, the density of the generalized Linnik distribution admits a
simple integral representation via known elementary densities (3),
(26) and (32).

By the way, it must be noted that any univariate symmetric random
variable $Y_{\alpha}$ is geometrically stable if and only if it is
representable as
$$
Y_{\alpha}=W_1^{1/\alpha}\circ S(\alpha, 0),\ \ \ 0<\alpha\le 2.
$$
Correspondingly, any univariate positive random variable
$Y_{\alpha}$ is geometrically stable if and only if it is
representable as
$$
Y_{\alpha}=W_1^{1/\alpha}\circ S(\alpha, 1),\ \ \ 0<\alpha\le 1.
$$
These representations immediately follow from the definition of
geometrically stable distributions and the transfer theorem for
cumulative geometric random sums, see \cite{GnedenkoKorolev1996}.

Hence, if $\nu\neq 1$, then from the identifiability of scale
mixtures of stable laws (see Lemma 1) it follows that the
generalized Linnik distribution and the generalized Mittag-Leffler
distributions are not geometrically stable.

Let $\Sigma$ be a positive definite $(r\times r)$-matrix,
$\alpha\in(0,2]$, $\nu>0$. As the `ordinary' multivariate Linnik
distribution, the multivariate generalized Linnik distribution can
be defined in at least two equivalent ways. First, it can be defined
by its characteristic function. Namely, a multivariate distribution
is called (centered elliptically contoured) generalized Linnik law,
if the corresponding characteristic function has the form
$$
\mathfrak{f}^{(L)}_{\alpha,\Sigma,\nu}(\emph{\textbf{t}})=\big[1+
(\emph{\textbf{t}}^{\top}\Sigma\emph{\textbf{t}})^{\alpha/2}\big]^{-\nu},\
\ \
\emph{\textbf{t}}\in\r^r.\eqno(38)   %(23)
$$
Second, let $\emph{\textbf{X}}$ be a random vector such that
$\mathcal{L}(\emph{\textbf{X}})=\mathfrak{N}_{\Sigma}$, independent
of the random variable $M_{\alpha/2,\nu}$ with the generalized
Mittag-Leffler distribution. By analogy with (33), introduce the
random vector $\emph{\textbf{L}}_{\alpha,\Sigma,\nu}$ as
$$
\emph{\textbf{L}}_{\alpha,\Sigma,\nu}=\sqrt{2M_{\alpha/2,\nu}}\circ\emph{\textbf{X}}.
$$
Then, in accordance with what has been said in Section 2.3,
$$
\mathcal{L}(\emph{\textbf{L}}_{\alpha,\Sigma,\nu})={\sf
E}\mathfrak{N}_{2M_{\alpha/2,\nu}\Sigma}.\eqno(39)   %(26)
$$
The distribution (29) will be called the $($centered$)$ elliptically
contoured multivariate generalized Linnik distribution.

Using Remark 1 we can easily make sure that the two definitions of
the multivariate generalized Linnik distribution coincide. Indeed,
with the account of (21), according to Remark 1, the characteristic
function of the random vector
$\emph{\textbf{L}}_{\alpha,\Sigma,\nu}$ defined by (39) has the form
$$
{\sf
E}\exp\{i\emph{\textbf{t}}^{\top}\emph{\textbf{L}}_{\alpha,\Sigma,\nu}\}=
\psi^{(M)}_{\alpha/2,\nu}\big(\emph{\textbf{t}}^{\top}\Sigma\emph{\textbf{t}}\big)=
\big[1+(\emph{\textbf{t}}^{\top}\Sigma\emph{\textbf{t}})^{\alpha/2}\big]^{-\nu}=
\mathfrak{f}^{(L)}_{\alpha,\Sigma,\nu}(\emph{\textbf{t}}),\ \
\emph{\textbf{t}}\in\r^r,
$$
that coincides with (38).

Based on (32), one more equivalent definition of the multivariate
Linnik distribution can be proposed. Namely, let
$\emph{\textbf{L}}_{\alpha,\Sigma,\nu}$ be an $r$-variate random
vector such that
$$
\emph{\textbf{L}}_{\alpha,\Sigma,\nu}=G_{\nu,1}^{1/\alpha}\circ\emph{\textbf{S}}(\alpha,\Sigma,0).\eqno(40)
$$
It is well known that the Laplace---Stieltjes transform
$\psi^{(G)}_{\nu,1}(s)$ of the random variable $G_{\nu,1}$ having
the gamma distribution with the shape parameter $\nu$ has the form
$$
\psi^{(G)}_{\nu,1}(s)=(1+s)^{-\nu},\ \ \ \ s>0.
$$
Then in accordance with Remark 1 the characteristic function of the
random vector $\emph{\textbf{L}}_{\alpha,\Sigma,\nu}$ defined by
(40) again has the form
$$
{\sf
E}\exp\{i\emph{\textbf{t}}^{\top}\emph{\textbf{L}}_{\alpha,\Sigma,\nu}\}=
\psi^{(G)}_{\nu,1}\big((\emph{\textbf{t}}^{\top}\Sigma\emph{\textbf{t}})^{\alpha/2}\big)=
\big[1+(\emph{\textbf{t}}^{\top}\Sigma\emph{\textbf{t}})^{\alpha/2}\big]^{-\nu}=
\mathfrak{f}^{(L)}_{\alpha,\Sigma,\nu}(\emph{\textbf{t}}),\ \
\emph{\textbf{t}}\in\r^r.
$$

Definitions (39) and (40) open the way to formulate limit theorems
stating that the multivariate generalized Linnik distribution can be
limiting both for random sums of independent identically distributed
random vectors with infinite second moments, and for random sums of
independent random vectors with finite covariance matrices.

\section{Scale-mixed stable distributions}

\subsection{Definition and general properties}

Let $\alpha\in(0,2]$, let $U$ be a positive random variable and
$\Sigma$ be a positive definite $(r\times r)$-matrix. An $r$-variate
random vector $\emph{\textbf{Y}}_{\alpha,\Sigma,0}$ is said to have
the {\it $U$-scale-mixed centered elliptically contoured stable
distribution}, if
$$
\mathcal{L}(\emph{\textbf{Y}}_{\alpha,\Sigma,0})={\sf
E}\mathfrak{S}_{\alpha,U^{2/\alpha}\Sigma,0}.
$$
In terms of random vectors this means that
$\emph{\textbf{Y}}_{\alpha,\Sigma,0}$ is representable as
$$
\emph{\textbf{Y}}_{\alpha,\Sigma,0}=U^{1/\alpha}\circ\emph{\textbf{S}}(\alpha,\Sigma,0).
$$
Correspondingly, for $0<\alpha\le1$, a univariate positive random
variable $Y_{\alpha,1}$ is said to to have the $U$-scale-mixed
one-sided stable distribution, if $Y_{\alpha,1}$ is representable as
$$
Y_{\alpha,1}=U^{1/\alpha}\circ S(\alpha,1).
$$

\smallskip

{\sc Theorem 1}. {\it Let $U$ be a positive random variable,
$\Sigma$ be a positive definite $(r\times r)$-matrix,
$\alpha\in(0,2]$, $\alpha'\in(0,1]$. Let}
$\emph{\textbf{S}}_{\alpha,\Sigma,0}$ {\it be a random vector with the
elliptically contoured symmetric stable distribution. Let an
r-variate random vector} $\emph{\textbf{Y}}_{\alpha\alpha',\Sigma,0}$
{\it have the $U$-scale-mixed symmetric stable distribution and a random
variable $Y_{\alpha',1}$ have the $U$-scale-mixed one-sided stable
distribution. Assume that} $\emph{\textbf{S}}_{\alpha,\Sigma,0}$ {\it and
$Y_{\alpha',1}$ are independent. Then}
$$
\emph{\textbf{Y}}_{\alpha\alpha',\Sigma,0}\eqd
Y_{\alpha',1}^{1/\alpha}\circ\emph{\textbf{S}}(\alpha,\Sigma,0).
$$

\smallskip

{\sc Proof}. From the definition of a $U$-scale-mixed stable
distribution and (11) we have
$$
\emph{\textbf{Y}}_{\alpha\alpha',\Sigma,0}\eqd
U^{1/\alpha\alpha'}\circ\emph{\textbf{S}}(\alpha\alpha',\Sigma,0)\eqd
U^{1/\alpha\alpha'}\circ
S^{1/\alpha}(\alpha',1)\circ\emph{\textbf{S}}(\alpha,\Sigma,0)\eqd
$$
$$
\eqd\big(U^{1/\alpha'}\circ
S(\alpha',1)\big)^{1/\alpha}\circ\emph{\textbf{S}}(\alpha,\Sigma,0)\eqd
Y_{\alpha',1}^{1/\alpha}\circ\emph{\textbf{S}}(\alpha,\Sigma,0).
$$

\smallskip

With $\alpha=2$, from Theorem 1 we obtain the following statement.

\smallskip

{\sc Corollary 1}. {\it Let $\alpha\in(0,2)$, $U$ be a positive
random variable, $\Sigma$ be a positive definite $(r\times
r)$-matrix, $\emph{\textbf{X}}$ be a random vector such that
$\mathcal{L}(X)=\mathfrak{N}_{\Sigma}$. Then}
$$
\emph{\textbf{Y}}_{\alpha,\Sigma,0}\eqd
\sqrt{2Y_{\alpha/2,1}}\circ\emph{\textbf{X}}.
$$

\smallskip

In other words, any multivariate scale-mixed symmetric stable
distribution is a scale mixture of multivariate normal laws. On the
other hand, since the normal distribution is stable with $\alpha=2$,
any multivariate normal scale mixture is a `trivial' multivariate
scale-mixed stable distribution.

\smallskip

To give particular examples of `non-trivial' scale-mixed stable
distributions, note that

\begin{itemize}
\item if $U\eqd W_1$, then $Y_{\alpha,1}\eqd M_{\alpha}$ and
$\emph{\textbf{Y}}_{\alpha,\Sigma,0}\eqd
\emph{\textbf{L}}_{\alpha,\Sigma}$;
\item if $U\eqd G_{\nu,1}$, then $Y_{\alpha,1}\eqd M_{\alpha,\nu}$
and $\emph{\textbf{Y}}_{\alpha,\Sigma,0}\eqd
\emph{\textbf{L}}_{\alpha,\Sigma,\nu}$;
\item if $U\eqd S(\alpha',1)$ with $0<\alpha'\le1$, then $Y_{\alpha,1}\eqd
S(\alpha\alpha',1)$ and $\emph{\textbf{Y}}_{\alpha,\Sigma,0}\eqd
\emph{\textbf{S}}(\alpha\alpha',\Sigma,0)$.

\end{itemize}

Among possible mixing distributions of the random variable $U$, we
will distinguish a special class that can play important role in
modeling observed regularities by heavy-tailed distributions.
Namely, assume that $V$ is a positive random variable and let
$$
U\eqd V\circ G_{\nu,1},
$$
that is, the distribution of $U$ is a scale mixture of gamma
distributions. We will denote the class of these distributions as
$\mathcal{G}^{(V)}$. This class is rather wide and besides the gamma
distribution and its particular cases (exponential, Erlang,
chi-square, etc.) with exponentially fast decreasing tail, contains,
for example, Pareto and Snedecor--Fisher laws with power-type
decreasing tail. In the last two cases the random variable $V$ is
assumed to have the corresponding gamma and inverse gamma
distributions, respectively.

For $\mathcal{L}(U)\in\mathcal{G}^{(V)}$ we have
$$
Y_{\alpha,1}\eqd (V\circ G_{\nu,1})^{1/\alpha}\circ S(\alpha,1)\eqd
V^{1/\alpha}\circ\big(G_{\nu,1}^{1/\alpha}\circ S(\alpha,1)\big)\eqd
V^{1/\alpha}\circ M_{\alpha,\nu}
$$
and
$$
\emph{\textbf{Y}}_{\alpha,\Sigma,0}\eqd (V\circ
G_{\nu,1})^{1/\alpha}\circ \emph{\textbf{S}}(\alpha,\Sigma,0)\eqd
V^{1/\alpha}\circ\big(G_{\nu,1}^{1/\alpha}\circ
\emph{\textbf{S}}(\alpha,\Sigma,0)\big)\eqd V^{1/\alpha}\circ
\emph{\textbf{L}}_{\alpha,\Sigma,\nu}.
$$
This means that with $\mathcal{L}(U)\in\mathcal{G}^{(V)}$, the
$U$-scale-mixed stable distributions are scale mixtures of the
generalized Mittag-Leffler and multivariate generalized Linnik laws.

Therefore, in what follows we will pay a special attention to
mixture representations of the generalized Mittag-Leffler and
multivariate generalized Linnik distributions. These representations
can be easily extended to any $U$-scale-mixed stable distributions
with $\mathcal{L}(U)\in\mathcal{G}^{(V)}$.

\subsection{Mixture representations for the multivariate generalized Linnik
distribution}

Mixture representations for the generalized Mittag-Leffler
distribution were considered in \cite{Korolevetal2019}. Some of them
were exposed in Section 2.4. Hence, we will focus on the
multivariate generalized Linnik distribution.

For $\alpha\in(0,2]$, $\alpha'\in(0,1)$ and $\nu>0$ using (40) and
(11) we obtain the following chain of relations:
$$
\emph{\textbf{L}}_{\alpha\alpha',\Sigma,\nu}\eqd
G_{\nu,1}^{1/\alpha\alpha'}\circ\emph{\textbf{S}}(\alpha\alpha',\Sigma,0)\eqd
S^{1/\alpha}(\alpha',1)\circ G_{\nu,1}^{1/\alpha\alpha'}\circ
\emph{\textbf{S}}(\alpha,\Sigma,0)\eqd
$$
$$
\eqd\big(S(\alpha',1)\circ\overline{G}_{\nu,\alpha',1}\big)^{1/\alpha}\circ\emph{\textbf{S}}(\alpha,\Sigma,0)\eqd
M_{\alpha',\nu}^{1/\alpha}\circ \emph{\textbf{S}}(\alpha,\Sigma,0).
$$
Hence, the following statement holds representing the generalized
Linnik distribution as a scale mixture of an elliptically contoured
multivariate symmetric stable law with any greater characteristic
parameter, the mixing distribution being the generalized
Mittag-Leffler law.

\smallskip

{\sc Theorem 2.} {\it If $\alpha\in(0,2]$, $\alpha'\in(0,1)$ and
$\nu>0$, then}
$$
\emph{\textbf{L}}_{\alpha\alpha',\Sigma,\nu}\eqd
M_{\alpha',\nu}^{1/\alpha}\circ
\emph{\textbf{S}}(\alpha,\Sigma,0).\eqno(41)
$$

\smallskip

Now let $\nu\in(0,1]$. From (32) and (4) it follows that
$$
\emph{\textbf{L}}_{\alpha,\Sigma,\nu}\eqd G_{\nu,1}^{1/\alpha}\circ
\emph{\textbf{S}}(\alpha,\Sigma,0)\eqd Z_{\nu,1}^{-1/\alpha}\circ
W_1^{1/\alpha}\circ\emph{\textbf{S}}(\alpha,\Sigma,0)\eqd
Z_{\nu,1}^{-1/\alpha}\eqd Z_{\nu,1}^{-1/\alpha}\circ
\emph{\textbf{L}}_{\alpha,\Sigma}
$$
yielding the following statement.

\smallskip

{\sc Theorem 3.} {\it If $\nu\in(0,1]$ and $\alpha\in(0,2]$, then}
$$
\emph{\textbf{L}}_{\alpha,\Sigma,\nu}\eqd Z_{\nu,1}^{-1/\alpha}\circ
\emph{\textbf{L}}_{\alpha,\Sigma}.\eqno(42)
$$

\smallskip

In other words, with $\nu\in(0,1]$ and $\alpha\in(0,2]$, the
multivariate generalized Linnik distribution is a scale mixture of
`ordinary' multivariate Linnik distributions.

Let $\alpha\in(0,2]$ and the random vector $\Lambda_{\Sigma}$ have
the multivariate Laplace distribution with some positive definite
$(r\times r)$-matrix $\Sigma$. In \cite{Korolev2016} it was shown
that if $\delta\in(0,1]$, then
$$
W_{\delta}\eqd W_1\circ S^{-1}(\delta,1).\eqno(43)
$$
Hence, it can be easily seen that
$$
\emph{\textbf{L}}_{\alpha,\Sigma}\eqd W_1^{1/\alpha}\circ
\emph{\textbf{S}}(\alpha,\Sigma,0)\eqd \sqrt{2W_{\alpha/2}\circ
S(\alpha/2,1)}\circ \emph{\textbf{X}}\eqd\sqrt{2W_1\circ
R_{\alpha/2}}\circ\emph{\textbf{X}}\eqd \sqrt{R_{\alpha/2}}\circ
\Lambda_{\Sigma}.\eqno(44)
$$
So, from Theorem 3 and (44) we obtain the following statement.

\smallskip

{\sc Corollary 2}. {\it If $\nu\in(0,1]$ and $\alpha\in(0,2]$, then
the multivariate generalized Linnik distribution is a scale mixture
of multivariate Laplace distributions}:
$$
\emph{\textbf{L}}_{\alpha,\Sigma,\nu}\eqd Z_{\nu,1}^{-1/\alpha}\circ
\sqrt{R_{\alpha/2}}\circ \Lambda_{\Sigma}.
$$

\smallskip

From (22) with $\nu=1$ and (43) it can be seen that
$$
\emph{\textbf{L}}_{\alpha,\Sigma}\eqd \sqrt{2M_{\alpha/2}}\circ
\emph{\textbf{X}}.
$$
Therefore we obtain one more corollary of Theorem 3 representing the
multivariate generalized Linnik distribution via `ordinary'
Mittag-Leffler distributions.

\smallskip

{\sc Corollary 3}. {\it If $\nu\in(0,1]$ and $\alpha\in(0,2]$, then}
$$
\emph{\textbf{L}}_{\alpha,\Sigma,\nu}\eqd Z_{\nu,1}^{-1/\alpha}\circ
\sqrt{2M_{\alpha/2}}\circ\emph{\textbf{X}}.
$$

\section{Convergence of the distributions of random sequences with
independent indices to multivariate scale-mixed stable
distributions}

\subsection{General transfer theorem for the distributions of
multivariate random sequences with independent random indices}

In applied probability it is a convention that a model distribution
can be regarded as well-justified or adequate, if it is an {\it
asymptotic approximation}, that is, if there exists a rather simple
limit setting (say, schemes of maximum or summation of random
variables) and the corresponding limit theorem in which the model
under consideration manifests itself as a limit distribution. The
existence of such limit setting can provide a better understanding
of real mechanisms that generate observed statistical regularities,
see e. g., \cite{GnedenkoKorolev1996}.

In this section we will prove some limit theorems presenting
necessary and sufficient conditions for the convergence of the
distributions of random sequences with independent random indices
(including sums of a random number of random vectors and
multivariate statistics constructed from samples with random sizes)
to scale mixtures of multivariate elliptically contoured stable
distributions. In the next section, as particular cases, conditions
will be obtained for the convergence of the distributions of random
sums of random vectors with both infinite and finite covariance
matrices to the multivariate generalized Linnik distribution.

Consider a sequence $\{\emph{\textbf{S}}_n\}_{n\ge1}$ of random
elements taking values in $\r^r$. Let $\Xi(\r^r)$ be the set of all
nonsingular linear operators acting from $\r^r$ to $\r^r$. The
identity operator acting from $\r^r$ to $\r^r$ will be denoted
$I_r$. Assume that there exist sequences $\{B_n\}_{n\ge1}$ of
operators from $\Xi(\r^r)$ and $\{\emph{\textbf{a}}_n\}_{n\ge1}$ of
elements from $\r^r$ such that
$$
\emph{\textbf{Y}}_n\equiv
B_n^{-1}(\emph{\textbf{S}}_n-\emph{\textbf{a}}_n)\Longrightarrow
\emph{\textbf{Y}}\ \ \ (n\to\infty)\eqno(45)
$$
where $\emph{\textbf{Y}}$ is a random element whose distribution
with respect to ${\sf P}$ will be denoted $H$, $H={\cal L}(Y)$.

Along with $\{\emph{\textbf{S}}_n\}_{n\ge1}$, consider a sequence of
integer-valued positive random variables $\{N_n\}_{n\ge1}$ such that
for each $n\ge1$ the random variable $N_n$ is independent of the
sequence $\{\emph{\textbf{S}}_k\}_{k\ge1}$. Let
$\emph{\textbf{c}}_n\in\r^r$, $D_n\in\Xi(\r^r)$, $n\ge1$. Now we
will formulate sufficient conditions for the weak convergence of the
distributions of the random elements
$\emph{\textbf{Z}}_n=D_n^{-1}(\emph{\textbf{S}}_{N_n}-\emph{\textbf{c}}_n)$
as $n\to\infty$.

For $\emph{\textbf{g}}\in\r^r$ denote
$\emph{\textbf{W}}_n(\emph{\textbf{g}})=D_n^{-1}(B_{N_n}\emph{\textbf{g}}+\emph{\textbf{a}}_{N_n}-\emph{\textbf{c}}_n)$.
By measurability of a random field we will mean its measurability as
a function of two variates, an elementary outcome and a parameter,
with respect to the Cartesian product of the $\sigma$-algebra ${\cal
A}$ and the Borel $\sigma$-algebra ${\cal B}(\r^r)$ of subsets of
$\r^r$.

In \cite{KorolevKossova1992, KorolevKossova1995} the following
theorem was proved which establishes sufficient conditions of the
weak convergence of multivariate random sequences with independent
random indices under operator normalization.

\smallskip

{\sc Theorem 4} \cite{KorolevKossova1992, KorolevKossova1995}. {\it
Let $\|D_n^{-1}\|\to\infty$ as $n\to\infty$ and let the sequence of
random variables $\{\|D_n^{-1}B_{N_n}\|\}_{n\ge1}$ be tight. Assume
that there exist a random element} $\emph{\textbf{Y}}$ {\it with
distribution $H$ and an $r$-dimensional random field}
$\emph{\textbf{W}}(\emph{\textbf{g}})$, $\emph{\textbf{g}}\in\r^r$,
{\it such that $(45)$ holds and}
$$
\emph{\textbf{W}}_n(\emph{\textbf{g}})\Longrightarrow
\emph{\textbf{W}}(\emph{\textbf{g}})\ \ \ \ \ (n\to\infty)
$$
{\it for $H$-almost all} $\emph{\textbf{g}}\in\r^r$. {\it Then the
random field} $\emph{\textbf{W}}(\emph{\textbf{g}})$ {\it is
measurable, linearly depends on} $\emph{\textbf{g}}$ {\it and}
$$
\emph{\textbf{Z}}_n \Longrightarrow
\emph{\textbf{W}}(\emph{\textbf{Y}})\ \ \ \ \ (n\to\infty),
$$
{\it where the random field} $\emph{\textbf{W}}(\cdot)$ {\it and the
random element} $\emph{\textbf{Y}}$ {\it are independent}.

\smallskip

Now consider a special case of the general limit setting and assume
that the normalization is scalar and the limit random vector
$\emph{\textbf{Y}}$ in (45) has a stable distribution. Namely, let
$\{b_n\}_{n\ge1}$ be an infinitely increasing sequence of positive
numbers and, instead of the general condition (45) assume that
$$
\mathcal{L}\big(b_n^{-1/\alpha}\emph{\textbf{S}}_n\big)\Longrightarrow
\mathfrak{S}_{\alpha,\Sigma,0}\eqno(46)
$$
as $n\to\infty$, where $\alpha\in(0,2]$ and $\Sigma$ is some
positive definite matrix. In other words, let
$$
b_n^{-1/\alpha}\emph{\textbf{S}}_n\Longrightarrow
\emph{\textbf{S}}(\alpha,\Sigma,0)\ \ \ \ (n\to\infty).
$$
Let $\{d_n\}_{n\ge1}$ be an infinitely increasing sequence of
positive numbers. As $\emph{\textbf{Z}}_n$ take the scalar
normalized random vector
$$\emph{\textbf{Z}}_n=d_n^{-1/\alpha}\emph{\textbf{S}}_{N_n}.$$
The following result can be considered as a generalization of the
main theorem of \cite{Korolev1997}.

\smallskip

{\sc Theorem 5}. {\it Let $N_n\to\infty$ in probability as
$n\to\infty$. Assume that the random vectors}
$\emph{\textbf{X}}_1,\emph{\textbf{X}}_2,\ldots$ {\it satisfy
condition $(46)$ with $\alpha\in(0,2]$ and a positive definite
matrix $\Sigma$. Then a distribution $F$ such that}
$$
{\cal L}(\emph{\textbf{Z}}_n)\Longrightarrow F\ \ \ (n\to\infty),
\eqno(47)
$$
{\it exists if and only if there exists a distribution function
$V(x)$ satisfying the conditions
\begin{enumerate}
\item[{\rm(i)}] $V(x)=0$ for $x<0$;
\item[{\rm(ii)}] for any $A\in{\cal B}(\r^r)$
$$
F(A)={\sf E}\mathfrak{S}_{\alpha,U^{2/\alpha}\Sigma,0}(A)=
\int_{0}^{\infty}\mathfrak{S}_{\alpha,u^{2/\alpha}\Sigma,0}(A)dV(u),\
\ x\in\r^1;
$$
\item[{\rm(iii)}] ${\sf P}(b_{N_n}<d_nx)\Longrightarrow V(x)$,
$n\to\infty$.
\end{enumerate}
}

\smallskip

{\sc Proof}. $\,$ {\it The `if' part}. We will essentially exploit
Theorem 4. For each $n\ge1$ set
$\emph{\textbf{a}}_n=\emph{\textbf{c}}_n=0$,
$B_n=D_n=d_n^{1/\alpha}I_r$. Let $U$ be a random variable with the
distribution function $V(x)$. Note that the conditions of the
theorem guarantee the tightness of the sequence of random variables
$$
\|D_n^{-1}B_{N_n}\|=(b_{N_n}/d_n)^{1/\alpha},\ \ n=1,2,\ldots
$$
implied by its weak convergence to the random variable
$U^{1/\alpha}$. Further, in the case under consideration we have
$\emph{\textbf{W}}_n(\emph{\textbf{g}})=(b_{N_n}/d_n)^{1/\alpha}\cdot
\emph{\textbf{g}}$, $\emph{\textbf{g}}\in\r^r$. Therefore, the
condition $N_n/d_n\Longrightarrow U$ implies
$\emph{\textbf{W}}_n(\emph{\textbf{g}})\Longrightarrow
U^{1/\alpha}\emph{\textbf{g}}$ for all $\emph{\textbf{g}}\in\r^r$.

Condition (46) means that in the case under consideration
$H=\mathfrak{S}_{\alpha,\Sigma,0}$. Hence, by Theorem 4
$\emph{\textbf{Z}}_n\Longrightarrow
U^{1/\alpha}\circ\emph{\textbf{S}}(\alpha,\Sigma,0)$ (recall that
the symbol $\circ$ stands for the product of {\it independent}
random elements). The distribution of the random element
$U^{1/\alpha}\circ\emph{\textbf{S}}(\alpha,\Sigma,0)$ coincides with
${\sf E}\mathfrak{S}_{\alpha,U^{2/\alpha}\Sigma,0}$, see Section
2.3.

\smallskip

{\it The `only if' part}. Let condition (47) hold. Make sure that
the sequence $\{\|D_n^{-1}B_{N_n}\|\}_{n\ge1}$ is tight. Let
$\emph{\textbf{Y}}\eqd \emph{\textbf{S}}(\alpha,\Sigma,0)$. There
exist $\delta>0$ and $\rho>0$ such that
$$
{\sf P}(\|\emph{\textbf{Y}}\|>\rho)>\delta.\eqno(48)
$$
For $\rho$ specified above and an arbitrary $x>0$ we have
$$
{\sf P}(\|\emph{\textbf{Z}}_n\|>x)\ge{\sf
P}\big(\big\|d_n^{-1/\alpha}\emph{\textbf{S}}_{N_n}\big\|>x;\
\big\|b_{N_n}^{-1/\alpha}\emph{\textbf{S}}_{N_n}\big\|>\rho\big)=
$$
$$
={\sf
P}\big((b_{N_n}/d_n)^{1/\alpha}>x\cdot\big\|b_{N_n}^{-1/\alpha}\emph{\textbf{S}}_{N_n}\big\|^{-1};\
\big\|b_{N_n}^{-1/\alpha}\emph{\textbf{S}}_{N_n}\big\|>\rho\big)\ge
$$
$$
\ge {\sf P}\big((b_{N_n}/d_n)^{1/\alpha}>x/\rho;\
\big\|b_{N_n}^{-1/\alpha}\emph{\textbf{S}}_{N_n}\big\|>\rho\big)=
$$
$$ =\sum\nolimits_{k=1}^{\infty}{\sf P}(N_n=k){\sf
P}\big((b_k/d_n)^{1/\alpha}>x/\rho;\
\big\|b_k^{-1/\alpha}\emph{\textbf{S}}_k\big\|>\rho\big)=
$$
$$
=\sum\nolimits_{k=1}^{\infty}{\sf P}(N_n=k){\sf
P}\big((b_k/d_n)^{1/\alpha}> x/\rho\big){\sf
P}\big(\big\|b_k^{-1/\alpha}\emph{\textbf{S}}_k\big\|>\rho\big)\eqno(49)
$$
(the last equality holds since any constant is independent of any
random variable). Since by (46) the convergence
$b_k^{-1/\alpha}\emph{\textbf{S}}_k\Longrightarrow Y$ takes place as
$k\to\infty$, from (48) it follows that there exists a number
$k_0=k_0(\rho,\delta)$ such that
$$
{\sf
P}\big(\big\|b_k^{-1/\alpha}\emph{\textbf{S}}_k\big\|>\rho\big)>\delta/2
$$
for all $k>k_0$. Therefore, continuing (49) we obtain
$$
{\sf
P}(\|\emph{\textbf{Z}}_n\|>x)\ge\frac{\delta}{2}\sum\nolimits_{k=k_0+1}^{\infty}
{\sf P}(N_n=k){\sf P}\big((b_k/d_n)^{1/\alpha}>x/\rho\big)=
$$
$$
=\frac{\delta}{2}\big[{\sf
P}\big((b_{N_n}/d_n)^{1/\alpha}>x/\rho\big)-
\sum\nolimits_{k=1}^{k_0}{\sf P}(N_n=k){\sf
P}\big((b_k/d_n)^{1/\alpha}>x/\rho\big)\big]\ge
$$
$$
\ge\frac{\delta}{2}\big[{\sf
P}\big((b_{N_n}/d_n)^{1/\alpha}>x/\rho\big)- {\sf P}(N_n\le
k_0)\big].
$$
Hence,
$$ {\sf
P}\big((b_{N_n}/d_n)^{1/\alpha}>x/R\big)\le\frac{2}{\delta} {\sf
P}(\|\emph{\textbf{Z}}_n\|>x)+{\sf P}(N_n\le k_0).\eqno(50)
$$
From the condition $N_n\pto\infty$ as $n\to\infty$ it follows that
for any $\epsilon>0$ there exists an $n_0=n_0(\epsilon)$ such that
${\sf P}(N_n\le n_0)<\epsilon$ for all $n\ge n_0$. Therefore, with
the account of the tightness of the sequence
$\{\emph{\textbf{Z}}_n\}_{n\ge1}$ that follows from its weak
convergence to the random element $\emph{\textbf{Z}}$ with
$\mathcal{L}(\emph{\textbf{Z}})=F$ implied by (47), relation (50)
implies
$$
\lim_{x\to\infty}\sup_{n\ge n_0(\epsilon)}{\sf
P}\big((b_{N_n}/d_n)^{1/\alpha}> x/\rho\big)\le\epsilon,\eqno(51)
$$
whatever $\epsilon>0$ is. Now assume that the sequence
$$
\|D_n^{-1}B_{N_n}\|=(b_{N_n}/d_n)^{1/\alpha},\ \ n=1,2,\ldots
$$
is not tight. In that case there exists an $\gamma>0$ and sequences
${\cal N}$ of natural and $\{x_n\}_{n\in{\cal N}}$ of real numbers
satisfying the conditions $x_n\uparrow\infty$ $(n\to\infty,\
n\in{\cal N})$ and
$$
{\sf P}\big((b_{N_n}/d_n)^{1/\alpha}>x_n\big)>\gamma,\ \ n\in{\cal
N}. \eqno(52)
$$
But, according to (51), for any $\epsilon>0$ there exist
$M=M(\epsilon)$ and $n_0=n_0(\epsilon)$ such that
$$
\sup_{n\ge n_0(\epsilon)}{\sf P}\big((b_{N_n}/d_n)^{1/\alpha}>
M(\epsilon)\big)\le 2\epsilon.\eqno(53)
$$
Choose $\epsilon<\gamma/2$ where $\gamma$ is the number from (52).
Then for all $n\in{\cal N}$ large enough, in accordance with (52),
the inequality opposite to (53) must hold. The obtained
contradiction by the Prokhorov theorem proves the tightness of the
sequence $\{\|D_n^{-1}B_{N_n}\|\}_{n\ge1}$ or, which in this case is
the same, of the sequence $\{b_{N_n}/d_n\}_{n\ge1}$.

Introduce the set ${\cal W}(\emph{\textbf{Z}})$ containing all
nonnegative random variables $U$ such that ${\sf
P}(\emph{\textbf{Z}}\in A)={\sf
E}\mathfrak{S}_{\alpha,U^{2/\alpha}\Sigma,0}(A)$ for any $A\in{\cal
B}(\r^r)$. Let $\lambda(\cdot,\cdot)$ be any probability metric that
metrizes weak convergence in the space of $r$-variate random
vectors, or, which is the same in this context, in the space of
distributions, say, the L\'evy--Prokhorov metric. If
$\emph{\textbf{X}}_1$ and $\emph{\textbf{X}}_2$ are random variables
with the distributions $F_1$ and $F_2$ respectively, then we
identify $\lambda(\emph{\textbf{X}}_1,\emph{\textbf{X}}_2)$ and
$\lambda(F_1,F_2)$). Show that there exists a sequence of random
variables $\{U_n\}_{n\ge1}$, $U_n\in{\cal W}(\emph{\textbf{Z}})$,
such that
$$
\lambda\big(b_{N_n}/d_n,\,U_n \big)\longrightarrow 0\ \
(n\to\infty).\eqno(54)
$$
Denote
$$
\beta_n=\inf\big\{\lambda\big(b_{N_n}/d_n,\,U\big):\ U\in{\cal
W}(\emph{\textbf{Z}})\big\}.
$$
Prove that $\beta_n\to 0$ as $n\to\infty$. Assume the contrary. In
that case $\beta_n\ge\delta$ for some $\delta>0$ and all $n$ from
some subsequence ${\cal N}$ of natural numbers. Choose a subsequence
${\cal N}_1\subseteq{\cal N}$ so that the sequence
$\{b_{N_n}/d_n\}_{n\in {\cal N}_1}$ weakly converges to a random
variable $U$ (this is possible due to the tightness of the family
$\{b_{N_n}/d_n\}_{n\ge1}$ established above). But then
$\emph{\textbf{W}}_n(\emph{\textbf{g}})\Longrightarrow
U^{1/\alpha}\emph{\textbf{g}}$ as $n\to\infty$, $n\in{\cal N}_1$ for
any $\emph{\textbf{g}}\in\r^r$. Applying Theorem 4 to $n\in{\cal
N}_1$ with condition (46) playing the role of condition (45), we
make sure that $U\in{\cal W}(\emph{\textbf{Z}})$, since condition
(47) provides the coincidence of the limits of all weakly convergent
subsequences. So, we arrive at the contradiction to the assumption
that $\beta_n\ge\delta$ for all $n\in{\cal N}_1$. Hence, $\beta_n\to
0$ as $n\to\infty$.

For any $n=1,2,\ldots$ choose a random variable $U_n$ from ${\cal
W}(\emph{\textbf{Z}})$ satisfying the condition
$$
\lambda\big(b_{N_n}/d_n,\,U_n\big)\le\beta_n+{\textstyle\frac{1}{n}}.
$$
This sequence obviously satisfies condition (54). Now consider the
structure of the set ${\cal W}(\emph{\textbf{Z}})$. This set
contains all the random variables defining the family of special
mixtures of multivariate centered elliptically contoured stable laws
considered in Lemma 1, according to which this family is
identifiable. So, whatever a random element $\emph{\textbf{Z}}$ is,
the set ${\cal W}(\emph{\textbf{Z}})$ contains at most one element.
Therefore, actually condition (54) is equivalent to
$$
b_{N_n}/d_n\Longrightarrow U\ \ \ (n\to\infty),
$$
that is, to condition (iii) of the theorem. The theorem is proved.

\smallskip

{\sc Corollary 4}. {\it Under the conditions of Theorem $5$,
non-randomly normalized random sequences with independent random
indices} $d_n^{-1/\alpha}\emph{\textbf{S}}_{N_n}$ {\it have the
limit stable distribution $\mathfrak{S}_{\alpha,\Sigma',0}$ with
some positive definite matrix $\Sigma'$ if and only if there exists
a number $c>0$ such that
$$
b_{N_n}/d_n\Longrightarrow c\ \ (n\to\infty).
$$
Moreover, in this case $\Sigma'=c^{2/\alpha}\Sigma$.}

\smallskip

This statement immediately follows from Theorem 5 with the account
of Lemma 1.

\subsection{Convergence of the distributions of random sums of random
vectors to special scale-mixed multivariate elliptically contoured
stable laws}

In Section 3 (see Corollary 1) we made sure that all scale-mixed
centered elliptically contoured stable distributions are
representable as multivariate normal scale mixtures. Together with
Theorem 5 this observation allows to suspect at least two
principally different limit schemes in which each of these
distributions can appear as limiting for random sums of independent
random vectors. We will illustrate these two cases by the example of
the multivariate generalized Linnik distribution.

As we have already mentioned, `ordinary' Linnik distributions are
geometrically stable. Geometrically stable distributions are only
possible limits for the distributions of geometric random sums of
independent identically distributed random vectors. As this is so,
the distributions of the summands belong to the domain of attraction
of the multivariate strictly stable law with some characteristic
exponent $\alpha\in(0,2]$ and hence, for $0<\alpha<2$ the univariate
marginals have infinite moments of orders greater or equal to
$\alpha$. As concerns the case $\alpha=2$, where the variances of
marginals are finite, within the framework of the scheme of
geometric summation in this case the only possible limit law is the
multivariate Laplace distribution \cite{Klebanov}.

Correspondinly, as we will demonstrate below, the multivariate
generalized Linnik distributions can be limiting for negative
binomial sums of independent identically distributed random vectors.
Negative binomial random sums turn out to be important and adequate
models of characteristics of precipitation (total precipitation
volume, etc.) during wet (rainy) periods in meteorology \cite{Gulev,
KorolevGorsheninDAN}. However, in this case the summands (daily
rainfall volumes) also must have distributions from the domain of
attraction of a strictly stable law with some characteristic
exponent $\alpha\in(0,2]$ and hence, with $\alpha\in(0,2)$, have
infinite variances, that seems doubtful, since to have an infinite
variance, the random variable {\it must} be allowed to take {\it
arbitrarily large} values with {\it positive} probabilities. If
$\alpha=2$, then the only possible limit distribution for negative
binomial random sums is the so-called variance gamma distribution
which is well known in financial mathematics
\cite{GnedenkoKorolev1996}.

However, when the (generalized) Linnik distributions are used as
models of statistical regularities observed in real practice and an
additive structure model is used of type of a (stopped) random walk
for the observed process, the researcher cannot avoid thinking over
the following question: which of the two combinations of conditions
can be encountered more often:

\begin{itemize}

\item the distribution of the number of summands (the number of jumps of a random walk)
is asymptotically gamma (say, negative binomial), but the
distributions of summands (jumps) have so heavy tails that, at
least, their variances are infinite, or

\item the second moments (variances) of the summands (jumps) are finite, but the
number of summands exposes an irregular behavior so that its very
large values are possible?

\end{itemize}

Since, as a rule, when real processes are modeled, there are no
serious reasons to reject the assumption that the variances of jumps
are finite, the second combination at least deserves a thorough
analysis.

As it was demonstrated in the preceding section, the scale-mixed
multivariate elliptically contoured stable distributions (including
multivariate (generalized) Linnik laws) even with $\alpha<2$ can be
represented as multivariate normal scale mixtures. This means that
they can be limit distributions in analogs of the central limit
theorem for random sums of independent random vectors {\it with
finite covariance matrices}. Such analogs with univariate `ordinary'
Linnik limit distributions were presented in
\cite{KorolevZeifmanKMJ} and extended to generalized Linnik
distributions in \cite{Korolevetal2019}. In what follows we will
present some examples of limit settings for random sums of
independent random vectors with principally different tail behavior.
In particular, it will de demonstrated that the scheme of negative
binomial summation is far not the only asymptotic setting (even for
sums of independent random variables!) in which the multivariate
generalized Linnik law appears as the limit distribution.

\smallskip

{\sc Remark 2}. Based on the results of
\cite{KorolevZeifman2016JSPI}, by an approach that slightly differs
from the one used here by the starting point, in the paper
\cite{Korchagin2015} it was demonstrated that if the random vectors
$\{\emph{\textbf{S}}_n\}_{n\ge1}$ are formed as cumulative sums of
independent random vectors:
$$
\emph{\textbf{S}}_n =\emph{\textbf{X}}_1+\ldots
+\emph{\textbf{X}}_n\eqno(55)
$$
for $n\in\mathbb{N}$, where $\emph{\textbf{X}}_1,
\emph{\textbf{X}}_2,\ldots$ are independent $r$-valued random
vectors, then the condition $N_n\pto\infty$ in the formulations of
Theorem 5 and Corollary 4 can be omitted.

\smallskip

Throughout this section we assume that the random vectors
$\emph{\textbf{S}}_n$ have the form (55).

Let $U\in\cal{U}$ (see Section 2.3), $\alpha\in(0,2]$, $\Sigma$ be a
positive definite matrix. In Section 3.1 the $r$-variate random
vector $\emph{\textbf{Y}}_{\alpha,\Sigma,0}$ with the the
multivariate $U$-scale-mixed centered elliptically contoured stable
distribution was introduced as
$\emph{\textbf{Y}}_{\alpha,\Sigma,0}=U^{1/\alpha}\circ\emph{\textbf{S}}(\alpha,\Sigma,0)$.
In this section we will consider the conditions under which
multivariate $U$-scale-mixed stable distributions can be limiting
for sums of independent random vectors.

%Again consider the sequence $\emph{\textbf{X}}_1,
%\emph{\textbf{X}}_2,\ldots$ of independent (in general, not
%necessarily identically distributed) $r$-valued random vectors, and
%for $n\in\mathbb{N}$ let $\emph{\textbf{S}}_n
%=\emph{\textbf{X}}_1+\ldots +\emph{\textbf{X}}_n$.
Consider a sequence of integer-valued positive random variables
$\{N_n\}_{n\ge1}$ such that for each $n\ge1$ the random variable
$N_n$ is independent of the sequence
$\{\emph{\textbf{S}}_k\}_{k\ge1}$. First, let $\{b_n\}_{n\ge1}$ be
an infinitely increasing sequence of positive numbers such that
convergence (46) takes place. Let $\{d_n\}_{n\ge1}$ be an infinitely
increasing sequence of positive numbers. The following statement
presents necessary and sufficient conditions for the convergence
$$
d_n^{-1/\alpha}\emph{\textbf{S}}_{N_n}\Longrightarrow
\emph{\textbf{Y}}_{\alpha,\Sigma,0}\ \ \ \ \ (n\to\infty).\eqno(56)
$$

\smallskip

{\sc Theorem 6.} {\it Under condition $(46)$, convergence $(56)$
takes place if and only if}
$$
b_{N_n}/d_n\Longrightarrow U\ \ \ \ \ (n\to\infty).\eqno(57)
$$

\smallskip

{\sc Proof}. This theorem is a direct consequence of Theorem 5 and
the definition of $\emph{\textbf{Y}}_{\alpha,\Sigma,0}$ with the
account of Remark 2.

\smallskip

{\sc Corollary 5}. {\it Assume that $\nu>0$. Under condition $(46)$,
the convergence}
$$
d_n^{-1/\alpha}\emph{\textbf{S}}_{N_n}\Longrightarrow
\emph{\textbf{L}}_{\alpha,\Sigma,\nu}\ \ \ \ \
(n\to\infty)%.\eqno(56)
$$
{\it takes place if and only if}
$$
b_{N_n}/d_n\Longrightarrow G_{\nu,1}\ \ \ \ \ (n\to\infty).\eqno(58)
$$

\smallskip

{\sc Proof}. To prove this statement it suffices to notice that the
multivariate generalized Linnik distribution is a $U$-scale-mixed
stable distribution with $U\eqd G_{\nu,1}$ (see representation (40))
and refer to Theorem 6 with the account of Remark 2.

\smallskip

Condition (58) holds, for example, if $b_n=d_n=n$, $n\in\mathbb{N}$,
and the random variable $N_n$ has the negative binomial distribution
with shape parameter $\nu>0$, that is, $N_n=\nb_{\nu,p_n}$,
$$
{\sf
P}(\nb_{\nu,p_n}=k)=\frac{\Gamma(\nu+k-1)}{(k-1)!\Gamma(r)}\cdot
p_n^{\nu}(1-p_n)^{k-1},\ \ \ \ k=1,2,...,
$$
with $p_n=n^{-1}$ (see, e. g., \cite{BeningKorolev2005,
SchluterTrede2016}). In this case ${\sf E}\nb_{\nu,p_n}=n\nu$.

\smallskip

Now consider the conditions providing the convergence in
distribution of scalar normalized random sums of independent random
vectors satisfying condition (46) with some $\alpha\in(0,2]$ and
$\Sigma$ to a random vector $\emph{\textbf{Y}}_{\beta,\Sigma,0}$
with the $U$-scale-mixed stable distribution ${\sf
E}\mathfrak{S}_{\beta,U^{2/\beta}\Sigma,0}$ with some
$\beta\in(0,\alpha)$. For convenience, let $\beta=\alpha\alpha'$
where $\alpha'\in(0,1)$.

Recall that in Section 3.1, for $\alpha'\in(0,1]$ the positive
random variable $Y_{\alpha',1}$ with the univariate one-sided
$U$-scale-mixed stable distribution was introduced as
$Y_{\alpha',1}\eqd U^{1/\alpha}\circ S(\alpha',1)$.

\smallskip

{\sc Theorem 7}. {\it Let $\alpha'\in(0,1]$. Under condition $(46)$,
the convergence}
$$
d_n^{-1/\alpha}\emph{\textbf{S}}_{N_n}\Longrightarrow
\emph{\textbf{Y}}_{\alpha\alpha',\Sigma,0}\ \ \ \ \
(n\to\infty)%.\eqno(56)
$$
{\it takes place if and only if}
$$
b_{N_n}/d_n\Longrightarrow Y_{\alpha',1}\ \ \ \ \ (n\to\infty).%\eqno(58)
$$

\smallskip

{\sc Proof}. This statement directly follows from Theorems 1 and 5
with the account of Remark 2.

\smallskip

{\sc Corollary 6}. {\it Let $\alpha'\in(0,1]$, $\nu>0$. Under
condition $(46)$, the convergence}
$$
d_n^{-1/\alpha}\emph{\textbf{S}}_{N_n}\Longrightarrow
\emph{\textbf{L}}_{\alpha\alpha',\Sigma,\nu}\ \ \ \ \
(n\to\infty)%.\eqno(56)
$$
{\it takes place if and only if}
$$
b_{N_n}/d_n\Longrightarrow M_{\alpha',\nu}\ \ \ \ \ (n\to\infty).%\eqno(58)
$$

\smallskip

{\sc Proof}. This statement directly follows from Theorems 2 (see
representation (41)) and 7 with the account of Remark 2.

\smallskip

From the case of heavy tails turn to the `light-tails' case where in
(46) $\alpha=2$. In other words, assume that the properties of the
summands $\emph{\textbf{X}}_j$ provide the asymptotic normality of
the sums $\emph{\textbf{S}}_n$. More precisely, assume that instead
of (46), the condition
$$
b_n^{-1/2}\emph{\textbf{S}}_n\Longrightarrow \emph{\textbf{X}}\ \ \
\ \ (n\to\infty)\eqno(59)
$$
holds. The following results show that even under condition (59),
heavy-tailed $U$-scale-mixed multivariate stable distributions can
be limiting for random sums.

\smallskip

{\sc Theorem 8.} {\it Under condition $(59)$, convergence $(56)$
takes place if and only if}
$$
b_{N_n}/d_n\Longrightarrow Y_{\alpha/2,1}\ \ \ \ \
(n\to\infty).%\eqno(57)
$$

\smallskip

{\sc Proof}. This theorem is a direct consequence of Theorem 5 and
Corollary 1, according to which
$\emph{\textbf{Y}}_{\alpha,\Sigma,0}\eqd \sqrt{2Y_{\alpha/2,1}}\circ
\emph{\textbf{X}}$ with the account of Remark 2.

\smallskip

{\sc Corollary 7}. {\it Assume that $N_n\to\infty$ in probability as
$n\to\infty$. Under condition $(59)$, non-randomly normalized random
sums} $d_n^{-1/2}\emph{\textbf{S}}_{N_n}$ {\it have the limit stable
distribution $\mathfrak{S}_{\alpha,\Sigma,0}$ if and only if}
$$
b_{N_n}/d_n\Longrightarrow 2S(\alpha/2,1)\ \ \ \ \ (n\to\infty).
$$

\smallskip

{\sc Proof}. This statement follows from Theorem 8 with the account
of the univariate version of (10) (see Theorem 3.3.1 in
\cite{Zolotarev1986}) and Remark 2.

\smallskip

{\sc Corollary 8}. {\it Assume that $N_n\to\infty$ in probability as
$n\to\infty$, $\nu>0$. Under condition $(59)$, the convergence}
$$
d_n^{-1/2}\emph{\textbf{S}}_{N_n}\Longrightarrow
\emph{\textbf{L}}_{\alpha,\Sigma,\nu}\ \ \ \ \
(n\to\infty)%.\eqno(56)
$$
{\it takes place if and only if}
$$
b_{N_n}/d_n\Longrightarrow 2M_{\alpha/2,\nu}\ \ \ \ \
(n\to\infty).%\eqno(58)
$$

\smallskip

{\sc Proof}. To prove this statement it suffices to notice that the
multivariate generalized Linnik distribution is a multivariate
normal scale mixture with the generalized Mittag-Leffler mixing
distribution (see definition (39)) and refer to Theorem 8 with the
account of Remark 2.

Another way to prove Corollary 8 is to deduce it from Corollary 6.

\smallskip

Product representations for limit distributions in these theorems
proved in the preceding sections allow to use other forms of the
conditions for the convergence of random sums of random vectors to
particular scale mixtures of multivariate stable laws.

\renewcommand{\refname}{References}

\end{document}